
\def\LaTeX{%
  \let\Begin\begin
  \let\End\end
  \let\salta\relax
  \let\finqui\relax
  \let\futuro\relax}

\def\UK{\def\our{our}\let\sz s}
\def\USA{\def\our{or}\let\sz z}

\UK 



\LaTeX

\USA


\salta

\documentclass[twoside,12pt]{article}
\setlength{\textheight}{24cm}
\setlength{\textwidth}{16cm}
\setlength{\oddsidemargin}{2mm}
\setlength{\evensidemargin}{2mm}
\setlength{\topmargin}{-15mm}
\parskip2mm


\usepackage[usenames,dvipsnames]{color}
\usepackage{amsmath}
\usepackage{amsthm}
\usepackage{amssymb,bbm}
\usepackage[mathcal]{euscript}

\usepackage{cite}
\usepackage{hyperref}
\usepackage{enumitem}

\usepackage[ulem=normalem,draft]{changes}
%
%

%
 
\definecolor{ciclamino}{rgb}{0.5,0,0.5}
\definecolor{blu}{rgb}{0,0,0.7}
\definecolor{rosso}{rgb}{0.85,0,0}

\def\juerg #1{{\color{Green}#1}}
\def\pier #1{{\color{blue}#1}} 
\def\an #1{{\color{rosso}#1}}
\def\last #1{{\color{magenta}#1}}
\def\gianni #1{{\color{red}#1}}
\def\pcol #1{{\color{ciclamino}#1}} 

\def\oldgianni #1{#1}
\def\gianni #1{#1}
\def\pier #1{#1}
\def\an #1{#1}
\def\last #1{#1}
\def\juerg #1{#1}
\def\pcol #1{#1}




\bibliographystyle{plain}


%
\newtheorem{theorem}{Theorem}[section]
\newtheorem{remark}[theorem]{Remark}

\newtheorem{lemma}[theorem]{Lemma}

\finqui

\def\Beq{\Begin{equation}}
\def\Eeq{\End{equation}}
\def\Bsist{\Begin{eqnarray}}
\def\Esist{\End{eqnarray}}

\def\Bthm{\Begin{theorem}}
\def\Ethm{\End{theorem}}
\def\Blem{\Begin{lemma}}
\def\Elem{\End{lemma}}

\def\Brem{\Begin{remark}\rm}
\def\Erem{\End{remark}}

\def\Bcenter{\Begin{center}}
\def\Ecenter{\End{center}}
\let\non\nonumber




\def\step #1 \par{\medskip\noindent{\bf #1.}\quad}
\def\jstep #1: \par {\vspace{2mm}\noindent\underline{\sc #1 :}\par\nobreak\vspace{1mm}\noindent}

\def\aand{\quad\hbox{and}\quad}
\def\Lip{Lip\-schitz}
\def\Holder{H\"older}

\def\lhs{left-hand side}
\def\rhs{right-hand side}




\def\multibold #1{\def\arg{#1}%
  \ifx\arg\pto \let\next\relax
  \else
  \def\next{\expandafter
    \def\csname #1#1\endcsname{{\boldsymbol #1}}%
    \multibold}%
  \fi \next}

\def\pto{.}

\def\multical #1{\def\arg{#1}%
  \ifx\arg\pto \let\next\relax
  \else
  \def\next{\expandafter
    \def\csname cal#1\endcsname{{\cal #1}}%
    \multical}%
  \fi \next}

\def\multigrass #1{\def\arg{#1}%
  \ifx\arg\pto \let\next\relax
  \else
  \def\next{\expandafter
    \def\csname grass#1\endcsname{{\mathbb #1}}%
    \multigrass}%
  \fi \next}


\def\multimathop #1 {\def\arg{#1}%
  \ifx\arg\pto \let\next\relax
  \else
  \def\next{\expandafter
    \def\csname #1\endcsname{\mathop{\rm #1}\nolimits}%
    \multimathop}%
  \fi \next}

\multibold
qweryuiopasdfghjklzxcvbnmQWERTYUIOPASDFGHJKLZXCVBNM.  

\multical
QWERTYUIOPASDFGHJKLZXCVBNM.

\multigrass
QWERTYUIOPASDFGHJKLZXCVBNM.

\multimathop
diag dist div dom mean meas sign supp .

\def\Span{\mathop{\rm span}\nolimits}


\def\accorpa #1#2{\eqref{#1}--\eqref{#2}}
\def\Accorpa #1#2 #3 {\gdef #1{\eqref{#2}--\eqref{#3}}%
  \wlog{}\wlog{\string #1 -> #2 - #3}\wlog{}}


\def\separa{\noalign{\allowbreak}}

\def\somma #1#2#3{\sum_{#1=#2}^{#3}}

\def\graffe #1{\mathopen\{#1\mathclose\}}

\def\<#1>{\mathopen\langle #1\mathclose\rangle}
\def\norma #1{\mathopen \| #1\mathclose \|}

\def\aeO{\checkmmode{a.e.\ in~$\Omega$}}
\def\aeQ{\checkmmode{a.e.\ in~$Q$}}
\def\aet{\checkmmode{a.e.\ in~$(0,T)$}}
\def\aat{\checkmmode{for a.a.\ $t\in(0,T)$}}

\let\wt\widetilde
\def\cpto{\,\cdot\,}

\def\Beta{\widehat\beta}
\def\Pi{\widehat\pi}
\def\Feps{F_\eps}
\def\Betaeps{\Beta_\eps}
\def\betaeps{\beta_\eps}
\def\geps{g_\eps}
\def\heps{h_\eps}
\def\vveps{\vv_\eps}
\def\peps{p_\eps}
\def\phieps{\phi_\eps}
\def\mueps{\mu_\eps}
\def\sigeps{\sigma_\eps}

\def\ej{e_j}
\def\ezj{\ej^0}
\def\ei{e_i}
\def\ezi{\ei^0}
\def\lambdaj{\lambda_j}
\def\lambdazj{\lambdaj^0}
\def\Vk{V_k}
\def\Vzk{\Vk^0}
\def\vvk{\vv_k}
\def\uuk{\uu_k}
\def\qqk{\qq_k}
\def\pk{p_k}
\def\phik{\phi_k}
\def\muk{\mu_k}
\def\sigk{\sigma_k}
\def\phikj{\phi_{kj}}
\def\mukj{\mu_{kj}}
\def\sigkj{\sigma_{kj}}
\def\hphik{\an{\hat {\boldsymbol \phi}}_k}
\def\hmuk{\an{\hat {\boldsymbol \mu}}_k}
\def\hsigk{\an{\hat {\boldsymbol \sigma}}_k}
\def\vk{v_k}
\def\zk{z_k}
\def\zetak{\zeta_k}
\def\badphik{\badphi_k}

\def\iot {\int_0^t}
\def\ioT {\int_0^T}
\def\intQt{\int_{Q_t}}
\def\intQ{\int_Q}
\def\iO{\int_\Omega}
\def\iG{\int_\Gamma}
\def\intS{\int_\Sigma}
\def\intSt{\int_{\Sigma_t}}

\def\dt{\partial_t}
\def\dn{\partial_{\nn}}
\def\suG{|_\Gamma}

\def\checkmmode #1{\relax\ifmmode\hbox{#1}\else{#1}\fi}


\let\erre\grassR
\let\enne\grassN
\def\errebar{(-\infty,+\infty]}




\def\genspazio #1#2#3#4#5{#1^{#2}(#5,#4;#3)}
\def\spazio #1#2#3{\genspazio {#1}{#2}{#3}T0}

\def\L {\spazio L}
\def\H {\spazio H}
\def\W {\spazio W}


\def\Lx #1{L^{#1}(\Omega)}
\def\Hx #1{H^{#1}(\Omega)}
\def\Wx #1{W^{#1}(\Omega)}
\def\LxG #1{L^{#1}(\Gamma)}
\def\HxG #1{H^{#1}(\Gamma)}

\def\LQ #1{L^{#1}(Q)}
\def\LS #1{L^{#1}(\Sigma)}
\def\CQ #1{C^{#1}(\overline Q)}

\def\Luno{\Lx 1}

\def\Hdue{\Hx 2}
\def\Hunoz{{H^1_0(\Omega)}}

\def\LdueG{\LxG 2}


\def\LLx #1{\LL^{#1}(\Omega)}
\def\HHx #1{\HH^{#1}(\Omega)}
\def\WWx #1{\WW^{#1}(\Omega)}
\def\HHxG #1{\HH^{#1}(\Gamma)}


\let\badeps\epsilon
\let\eps\varepsilon
\let\badphi\phi
\let\phi\varphi

\let\sig\sigma

\let\TeXchi\chi                         
\newbox\chibox
\setbox0 \hbox{\mathsurround0pt $\TeXchi$}
\setbox\chibox \hbox{\raise\dp0 \box 0 }
\def\chi{\copy\chibox}



\def\normaV #1{\norma{#1}_V}

\def\phiz{\phi_0}
\def\sigz{\sigma_0}

\let\hat\widehat

\def\0{{\boldsymbol {0} }}
\let\T\grassT
\let\I\grassI
\def\Sph {{\cal S}_{\phi}}
\def\Ss {{\cal S}_{\sig}}
\def\Nph{N_\phi}
\def\Ns{N_\sig} 

\def\bzeta{{\boldsymbol \zeta}}
\def\h{\mathbbm{h}}

\let\emb\hookrightarrow

\def\Vz{V_0}
\def\Vp{{V^*}}

\def\Vzp{\Vz^*}

\def\soluz{(\vv,p,\mu,\phi,\xi,\sigma)}
\def\soluzeps{(\vveps,\peps,\mueps,\phieps,\sigeps)}
\def\soluzk{(\vvk,\pk,\muk,\phik,\sigk)}


\def\CP0{(${\mathcal{CP}}_0$)}

\usepackage{amsmath}
\DeclareFontFamily{U}{mathc}{}
\DeclareFontShape{U}{mathc}{m}{it}%
{<->s*[1.03] mathc10}{}

\DeclareMathAlphabet{\mathscr}{U}{mathc}{m}{it}

\Begin{document}


%
\title{Cahn--Hilliard--Brinkman model for tumor growth with possibly singular potentials}
\author{}
\date{}
\maketitle
\Bcenter
\vskip-1.5cm
{\large\sc Pierluigi Colli$^{(1,2)}$}\\
{\normalsize e-mail: {\tt pierluigi.colli@unipv.it}}\\[0.25cm]
{\large\sc Gianni Gilardi $^{(1,2)}$}\\
{\normalsize e-mail: {\tt gianni.gilardi@unipv.it}}\\[0.25cm]
{\large\sc Andrea Signori$^{(1)}$}\\
{\normalsize e-mail: {\tt andrea.signori01@unipv.it}}\\[0.25cm]
{\large\sc J\"urgen Sprekels$^{(3)}$}\\
{\normalsize e-mail: {\tt juergen.sprekels@wias-berlin.de}}\\[.5cm]
$^{(1)}$
{\small Dipartimento di Matematica ``F. Casorati''}\\
{\small Universit\`a di Pavia}\\
{\small via Ferrata 5, I-27100 Pavia, Italy}\\[.3cm] 
$^{(2)}$
{\small Research Associate at the IMATI -- C.N.R. Pavia}\\
{\small via Ferrata 5, I-27100 Pavia, Italy}\\[.3cm] 
$^{(3)}$
{\small Department of Mathematics}\\
{\small Humboldt-Universit\"at zu Berlin}\\
{\small Unter den Linden 6, D-10099 Berlin, Germany}\\[2mm]
{\small and}\\[2mm]
{\small Weierstrass Institute for Applied Analysis and Stochastics}\\
{\small Mohrenstrasse 39, D-10117 Berlin, Germany}\\[10mm]

\Ecenter
\Begin{abstract}
\noindent 
We analyze a \juerg{phase} field model for tumor growth consisting of a Cahn--Hilliard--Brinkman system, ruling the evolution of the tumor mass, coupled 
with an advection-reaction-diffusion equation for a chemical species acting as a nutrient.
The main novelty of the paper concerns the discussion of the existence 
of weak solutions to the system covering all the meaningful cases for the nonlinear potentials; in particular,
the typical choices given by the regular, the logarithmic, and the double obstacle potentials are 
\juerg{admitted in our treatise}.
Compared to previous results related to similar models, we suggest, instead of the classical no-flux condition, a
Dirichlet boundary condition for the chemical potential appearing in the Cahn--Hilliard-type equation.
Besides, abstract growth conditions for the source terms that may depend on the solution variables are postulated.


\End{abstract}
\salta
\pagestyle{myheadings}
\newcommand\testopari{\sc Colli -- Gilardi -- Signori -- Sprekels}
\newcommand\testodispari{\sc Cahn--Hilliard--Brinkman tumor model}
\markboth{\testopari}{\testodispari}
\finqui
%
\section{Introduction}
\label{INTRO}
\setcounter{equation}{0}
Cancer is \juerg{nowadays still} one of the \juerg{main diseases} causing death worldwide. \juerg{Beyond doubt,
the understanding of the development of solid tumor growth is one of the major} challenges scientists have to face in the current century.
Moreover, it is now, more than ever, apparent that 
\juerg{only interdisciplinary efforts may enable us to gain deeper insights into cancer development mechanisms.}
In this scenario, mathematics could play a crucial role, since multiscale mathematical modeling provides a quantitative
tool that may help in diagnostic and prognostic applications \juerg{(see, e.g., the seminal book \cite{CL})}. Among others, mathematics has two decisive advantages: the first one is that of being able to select particular mechanisms 
\juerg{that} may be more relevant than others, while the second one is that of being able to foresee and make predictions that may be precious for medical practitioners, without inflicting any harm to the patients.
Furthermore, the \juerg{extremely fast development} of computational methods \juerg{for the solution of} nonlinear PDEs 
opens the doors for a \juerg{direct} interaction between the experimental methods used by \juerg{physicians}
and the more theoretical \juerg{mathematical ones: indeed, advanced numerical} solvers may be implemented as a supporting tool in clinical therapies.
Recently, lots of phase field models modeling tumor growth have been proposed: a brief description of the state of the art
 will be provided later on. As biological materials like tumor agglomerates exhibit viscoelastic properties, we prescribe a velocity equation
of Brinkman type.

Let \gianni{$\Omega \subset \erre^3$} \juerg{be a spatial domain in which the tumor is located}, 
$T>0$ be a fixed final time, and set
$Q:= \Omega \times (0,T),$ and $\Sigma:= \partial \Omega \times (0,T)$.
Then the system \juerg{under investigation} in this paper 
is a Cahn--Hilliard--Brinkman model related to tumor growth and reads as follows:
\begin{alignat}{2}
	\label{eq:1}
	& - \div \gianni{\T(\vv,p)} + \nu \vv
	= \mu \nabla\phi + \big(\sig+ \chi(1-\phi)\big)\nabla\sig \qquad && \text{in $Q$,} 
	\\	\label{eq:2}
	& \div\vv  = g \qquad && \text{in $Q$,} 
	\\
	\label{eq:3}
	& {\dt \phi + \div(\phi \vv)  - \div(m(\phi) \nabla\mu)=  \Sph(\phi,\sig)} \qquad && {\text{in $Q$,}  }
	\\
	\label{eq:4}
	& \mu  = - \badeps \Delta \phi + \badeps^{-1 }F'(\phi) -\chi \sig \qquad && \text{in $Q$,} 
	\\
	\label{eq:5}
	& \dt \sig + \div\big(\sig \vv) - \div (n(\phi) \nabla(\sig+ \chi(1-\phi))\big) = \Ss(\phi,\sig) \qquad && \text{in $Q$,} 
\end{alignat}
\oldgianni{where the \juerg{coefficients} $m$ and $n$ are positive functions and} the \emph{viscous Brinkman stress tensor} 
$\T$ is defined~as
\begin{align}
	\gianni{\T(\vv,p)}
	& = 2 \eta D \vv
	+ \lambda \oldgianni{(\div\vv)} \I 
	- p \I \,.
	\label{def:T}
\end{align}
\oldgianni{Here, the standard notation 
\Beq
  D\vv := \frac 12 \bigl( \nabla\vv + (\nabla\vv)^T \bigr)
  \label{symgrad}
\Eeq
\Accorpa\Stress def:T symgrad
is used for the symmetrized gradient of the velocity field~$\vv$ which represents the \emph{volume-averaged velocity field} of the mixture with {\it permeability} $\nu$.
Moreover, $p$~is the pressure, $\eta$ and $\lambda$ are nonnegative \gianni{constants} denoting the
{\it shear viscosity} and the {\it bulk viscosity}, respectively,
and \juerg{$\I \in \erre^{3\times 3}$ is the identity matrix. Although}}
several choices are possible, we endow the above system with the following boundary and initial conditions:
\begin{alignat}{2}
	\label{eq:6}
	& \gianni{\T(\vv,p)}{\nn} = \0 \qquad && \text{on $\Sigma$,} 
	\\	
	\label{eq:7}
	& \dn \phi = 0,
	 \quad 
	 \mu =\mu_\Sigma , \quad
	\dn \sig =\kappa(\sig_\Sigma - \sigma) \qquad && \text{on $\Sigma$,} 
	\\
	\label{eq:8}
	& \phi(0) = \phiz, \quad \sig(0) = \sigz \quad  \qquad && \text{in $\Omega$,} 	
\end{alignat}
\Accorpa\Sys {eq:1} {eq:8}
where $\nn$ denotes the outer unit normal vector to $\partial \Omega$, and $\dn$ the associated outward normal derivative.
The other variables of the system are $\phi,\mu,$ and $ \sig.$
In relevant cases, the phase variable $\phi$ is an order parameter
taking values in $[-1,1]$ \juerg{that represents the difference between the volume fractions of tumor cells
and healthy cells}.
It allows us to keep track of the evolution of the \juerg{boundary} of the tumor, since
the level sets $\{\phi=1\}:=\{x \in \Omega: \phi(x)=1\}$ and $\{\phi=-1\}$ 
describe the region of pure phases: the tumorous phase and the healthy phase, respectively.
The second variable $\mu$ denotes the chemical potential related to $\phi$ as 
in the framework of the Cahn--Hilliard equation.
We postulate the growth and proliferation of the tumor to be driven by \juerg{the} absorption and consumption 
of some nutrient~\oldgianni{$\sig$} (usually oxygen).
\oldgianni{Whenever $0\leq\sig\leq1$}, $\sigma \simeq 1 $ represents a rich nutrient concentration,
whereas $\sigma \simeq 0$ a poor one.
The functions $m(\phi)$ and $n(\phi)$ are \juerg{nonnegative} mobility functions related to the phase and to the nutrient variables, respectively.
The \juerg{positive physical} constants $\badeps$ and $\nu$ are \pcol{connected} to the interfacial thickness and surface tension, while the 
nonnegative constant $\chi$ \pcol{stands for} the chemotactic sensitivity.
Finally, \gianni{$\Sph$ and $\Ss$} denote nonlinearities representing some source terms 
\juerg{that account for} the mutual interplay between tumor, healthy cells, and nutrients.
For \juerg{further} details concerning the modeling, we refer to \cite{GLSS} and the references therein.

Concerning the boundary conditions, we point out that \eqref{eq:6} can be understood as a {\it no-friction} boundary condition for the velocity field $\vv$.
This is a common request for similar systems, see, e.g., \cite{EGAR, EGAR2, KS2}, 
as it does not enforce any compatibility condition on the velocity source term
$g$ in \eqref{eq:2}\juerg{, which} is otherwise needed if one assumes \juerg{a} no-slip boundary condition 
like $\vv = \0$ on $\Sigma$ 
or \juerg{the} no-penetration boundary condition $\vv \cdot \nn =\0 $ on~$\Sigma$. 
Both these conditions entail that $\iO g=0$, which is not ideal from the modeling perspective.
Let us also notice that the vector field $\vv$ is not solenoidal 
\juerg{(as typically in fluid-type problems)}, 
\juerg{which entails challenges} from the mathematical viewpoint.

As a common denominator of \oldgianni{a more general} Cahn--Hilliard equation, 
in \eqref{eq:4} $F'$ represents the (generalized) derivative of a double-well shaped nonlinearity~$F$. 
Prototypical examples are the regular, the logarithmic,
and the double obstacle potentials. 
\juerg{These} read, in the order, as
\begin{align}
  &F_{reg}(r):=\frac14(r^2-1)^2\,, \quad r\in\erre\,,
  \label{Freg}
  \\
  \separa
  &F_{log}(r):=
  \begin{cases}
  		\frac{\theta}{2}\left[(1+r)\ln(1+r)+(1-r)\ln(1-r)\right]-\frac{\theta_0}{2}r^2
  		 & \text{if } r\in(-1,1)\,,
			\\
			\juerg{\theta\,\ln(2)-\frac{\theta_0}2}&  \juerg{       \text{if } r\in \{-1,1\}}\,,
  		\\		
  		+\infty &\text{otherwise}\,,
  \end{cases}
  \label{Flog}
  \\
  \separa
  &F_{dob}(r):=\begin{cases}
  c(1-r^2) \quad&\text{if } r\in[-1,1]\,,\\
  +\infty \quad&\text{otherwise}\,,
  \end{cases}
  \label{F2ob}
\end{align}%
for some positive \pcol{constants $\theta<\theta_0$ and $c$.}
Besides, in the case of 
nonregular potentials like the double obstacle \eqref{F2ob}, the second equation \eqref{eq:4}
has to be \gianni{read as a differential inclusion.}

\juerg{The main novelty of this paper} is the prescription of \juerg{the Dirichlet boundary condition 
$ \mu =\mu_\Sigma$ for the chemical potential  on $\Sigma$},
in contrast to the standard homogeneous Neumann (no-flux) boundary condition $\dn \mu =0$ on $\Sigma$ (see, e.g., \cite{EGAR,EGAR2,KS2}). 
The limitation behind the latter choice regards the 
nonlinear potentials that can be considered: in all \juerg{of} the aforementioned papers, the authors 
were forced to restrict the analysis \juerg{to regular potentials}, 
possibly of just quadratic growth at infinity\juerg{; singular potentials like \eqref{Flog} or \eqref{F2ob} were
excluded from the analysis. 
These, however, are} actually more relevant physically, since, if solutions to the system exist, 
then the condition $\phi \in [-1,1]$ is \juerg{automatically} fulfilled.
The restriction \juerg{of the admitted potentials originates from} the presence of the source term $\Sph$ in the 
Cahn--Hilliard equation \eqref{eq:1}.
Roughly speaking, \juerg{for proving} that $\mu \in \L2 {\Hx1}$, \juerg{the energetic approach 
provides just a control on $\nabla\mu$ in $L^2(Q)$}. 
The classical approach then requires the employment of the Poincar\'e--Wirtinger 
inequality along with a control of the spatial mean of $\mu$ in $L^2(0,T)$.
This can be \juerg{achieved} by comparison in equation \eqref{eq:4}, 
provided that $F$ is regular and its derivative $F'$ possesses a prescribed growth, 
\juerg{which leads to the choice of potentials of polynomial type}.
Therefore, the novelty of our work revolves around the different boundary condition,
which  allows us to apply another version of Poincar\'e's inequality to establish immediate control 
of $\mu \in \L2 {\Hx1}$ from the bound $\nabla\mu \in L^2(Q)$ without the need of an additional 
control of the spatial average of $\mu$\juerg{; in this way, unpleasant growth restrictions 
for the potential can be avoided} (see also \cite{GARL_2}).

Without claiming to be exhaustive, let us now review some literature connected with the system \Sys.
It is \juerg{well known} that the Brinkman law interpolates between the Stokes and \juerg{Darcy paradigms},
and it has become rather popular in recent times, see \cite{EGAR,EGAR2,KS2}.
For tumor growth models both descriptions seem reasonable, because the associated Reynolds number is very small. 
Formally, we recover the {\it Darcy limit} when $\eta\equiv\lambda\equiv0$, and $\nu >0$, where the boundary condition \eqref{eq:6} yields \juerg{that} $p=0$ on $\Sigma$; \juerg{similarly}, the {\it Stokes} limit is obtained when $\eta,\,\lambda>0$, and $\nu =0$. 
The Stokes equation was suggested, e.g., in \cite{Byr,Fri}, by approximating the tumor as
a viscous fluid, while \pcol{Darcy's} law describes a viscous fluid permeating a porous medium represented 
by the \juerg{extracellular matrix}
and \juerg{accounts} for the inclination of cells to move away from regions of high compression,
see, e.g., \cite{GLSS, GARL_4}. 
Stationary approximations for system \Sys\ are popular as well, and we mention \cite{GL1, GARL_2, FGR} and the references therein, where just polynomial-type potentials were considered. 
To cope with the case of singular potential, some authors \juerg{(see \cite{CGH})} suggested to include suitable relaxations.
Besides, we  refer to \cite{SS, FLR, FLS} and the references therein, for related {\it nonlocal} versions,  
to \juerg{\cite{GLS, ACHE,KRS}} for the additional coupling with elasticity, and to \cite{RSchS} for the 
coupling with the Keller--Segel equation. 

\subsection{Biological Examples and Modeling Considerations}
Before diving into the mathematical details, let us outline some physically meaningful choices for the source terms $\Sph$ and $\Ss$ introduced above. 

\noindent
{\bf Linear \pcol{k}inetic.}
A first typical form for $\Sph$ and $\Ss$ was proposed by Cristini et al. in \cite{CL} motivated by linear kinetic:
\begin{align}
	\label{lin:kin}
	\Sph(\phi,\sig) & :=  ({\calP}\sig - {\calA}) \h(\phi),
	\quad
	\Ss(\phi,\sig) := \calB (\sig_c - \sig) - \calC \sig \h(\phi).
\end{align}
In the above expressions, ${\calP}, \calA, \calB,$ and $\calC$ denote nonnegative constants 
related to biological quantities that are, in the order, the proliferation rate of tumoral cells by consumption of nutrient, the apoptosis rate, the consumption rate of the nutrient with respect to a
preexisting concentration $\sig_c$, and the nutrient consumption rate.
As for the function $\h$, it  denotes an interpolation function between $-1$ and $1$ with the property that $\h(-1)=0$
and $\h(1)=1$. Roughly speaking, $\h$ weights the corresponding mechanism compared to the amount of cancer located in that region
and ``turns off" the associated mechanism when the tumor is not present.

%

\noindent
{\bf Linear phenomenological laws for chemical reactions.} 
Another approach was proposed by accounting for linear phenomenological laws for chemical reactions  by
A. Hawkins-Daarud et al. in \cite{HD}, where the following form was suggested:
\begin{align}
	\label{phen:laws}
	\Sph(\phi,\sig,\mu)  =- \Ss(\phi,\sig,\mu) := P(\phi)(\sig + \chi(1-\phi) - \mu),
\end{align}
where $P$ stands for a suitable nonnegative proliferation function.

In the forthcoming analysis, we will proceed in an abstract fashion without prescribing explicit structures for the source terms, 
but just postulating suitable growth conditions.
Those are straightforwardly fulfilled by the above cases but this last one. Namely, 
we cannot allow a linear growth of the sources with respect to~$\mu$.
Nevertheless, let us claim that the special form \eqref{phen:laws} can still be considered
provided to adjust some estimates accordingly (cf. \cite{CGH, CGRS}). 

\bigskip
The above system \Sys\ (see, e.g., \cite{GLSS})  is naturally associated 
to a free energy $\cal E$ with following the structure:
\Beq
	\label{energy:GL}
	{\cal E}(\phi,\sig)  = \frac \badeps 2\iO |\nabla\phi|^2 + \frac 1\badeps\iO F(\phi) + \iO N(\phi,\sig),
\Eeq
\pcol{with $N$ specified by}
\Beq
	\label{defN}
	N(\phi,\sig)  =  \frac 12 |\sig|^2+ \chi \sig (1-\phi).
\Eeq
\pcol{Here,} the first two terms in ${\cal E}$ yield the well-known Ginzburg--Landau energy modeling phase segregation and adhesion effects, while 
$N$ stands for the chemical free energy density, respectively.
For convenience, let us immediately set a specific notation to denote the partial derivatives of the free energy $N$
with respect to the variables $\phi$ and~$\sig$:
\begin{align}
	\Nph(\phi,\sig)  := \partial_\phi N(\phi,\sig) 
	= - \chi \sig
	\aand
	\Ns(\phi,\sig) := \partial_\sig N(\phi,\sig)
	= \sig + \chi (1-\phi),
	\label{idnut}
\end{align}
\Accorpa\Nutrient defN idnut
as a direct computation shows. 
This notation will turn very convenient for the estimation procedure later on.

As they will play any role from the mathematical viewpoint, we will set $m(\cdot)\equiv n(\cdot)\equiv 1$, and $\badeps=1$
in our discussion.
However, let us claim that our well-posedness result might be proven also for suitable nonconstant, 
albeit non-degenerate, mobilities (see \cite{EGAR}).


\section{Notation, Assumptions and Main Results}
\label{SEC:NOT}
\setcounter{equation}{0}

\juerg{Throughout the paper, $\Omega$ is a bounded and connected open subset of~$\erre^3$
(the two-dimensional case can be treated in the same way)
having a smooth} boundary $\Gamma:=\partial\Omega$.
In the following, $|\Omega|$~denotes the Lebesgue measure of~$\Omega$.
Similarly, we write $|\Gamma|$ for the two-dimensional Hausdorff measure of~$\Gamma$.
Given a final time $T>0$, we set, for every $t\in(0,T]$,
\begin{align}
	Q_t:=\Omega\times(0,t),
	\quad 
	\Sigma_t:=\Gamma\times(0,t),
	\quad 
	Q:=Q_T,
	\quad 
	\text{\last{and}}
	\quad 
	\Sigma:=\Sigma_T.
	\label{defQt-etc}
\end{align}
Given a Banach space $X$, we denote by 
\oldgianni{$\norma\cpto_X$, $X^*$, and $\<{\cdot},{\cdot}>_X$,
its norm}, its dual space, and the associated duality pairing, respectively.
\oldgianni{As for the notation of \an{norms}, some exceptions \juerg{will be utilized} in the sequel.
Moreover, the symbol used for the norm in some space $X$ is adopted for the one in any power of~$X$.}
Another standard notation that we employ concerns vector\an{s}, or vector-valued functions, which are denoted by bold symbols; 
for instance, \an{$\0$} stands for the zero vector in $\erre^3$, and $\vv=(v_1,v_2,v_3)\in \erre^3$.
For $1\leq q\leq\infty$ and $k\geq 0$, we indicate the usual Lebesgue and 
Sobolev spaces on $\Omega$ by $L^q(\Omega)$ and $W^{k,q}(\Omega)$,
with the standard abbreviation $\Hx k:=\Wx{k,2}$.
The norm in $\Lx q$ is simply denoted by $\norma\cpto_q$,
and the same symbol is used for the norms in the analogous spaces constructed on $Q$, $\Gamma$ and $\Sigma$, 
if no confusion can arise.
Then, we introduce the shorthands
\gianni{%
\Bsist
  && H:=L^2(\Omega), \quad
  V:=H^1(\Omega), \quad
  \Vz := \Hunoz\an{,} 
  \non
  \\
  && \aand
  W:=\graffe{ v\in H^2(\Omega): \partial_{\nn}v=0\;\,\text{a.e.~on } \Gamma},
  \non
\Esist
and endow \juerg{these spaces} with their natural norms.
For simplicity, we write $\norma\cpto$ instead of~$\norma\cpto_H$.}
Besides, the space $H$ will be identified \juerg{with} its dual, so that 
we have the following continuous, dense, and compact embeddings:
\Beq
  W \emb V \emb H \emb V^* ,
  \non
\Eeq
\gianni{yielding that $(V,H,\Vp)$ is a Hilbert triplet.
Similarly, $(\Vz,H,\Vzp)$ is a Hilbert triplet \juerg{that will be used as well}.}
\gianni{We observe at once the compatible embeddings
\Bsist
  && V \emb \Hx{s_2} \emb \Hx{s_1} \emb H
  \emb (\Hx{s_1})^* \emb (\Hx{s_2})^* \emb \Vp
  \non
  \\
  && \quad \hbox{for $0<s_1<s_2<1$}, 
  \label{embeddings}
\Esist
and we recall that
\Beq
  (\Hx s)^*
  = (H^s_0(\Omega))^*
  = \Hx{-s}
  \quad \hbox{for $0\leq s\leq 1/2$} \,.
  \label{dualeHs}
\Eeq
}%
\oldgianni{Finally, the same symbols written with boldface characters 
denote the corresponding spaces of vector\an{-}valued functions.
So, we have for instance that
$\LL^p(\Omega)=(\Lx p)^3$, $\HH=H^3$ and $\VV=V^3$.}
\gianni{However, since $\HH$ and $\VV$ are powers of $H$ and~$V$, 
we simply write $\norma\cpto$ and $\normaV\cpto$ 
instead of $\norma\cpto_{\HH}$ and $\norma\cpto_{\VV}$, respectively.}
\an{Moreover, for given matrices ${\AA,\BB}\in \erre^{3\times 3 }$, we define the scalar product
\begin{align*}
	{\AA : \BB} \,:= \sum_{i=1}^3 \sum_{j=1}^3 [{\AA}]_{ij}\, [{\BB}]_{ij} \;.
\end{align*}}

The following structural assumptions will be in order in our analysis.
The double\an{-}well potential $F$ \juerg{introduced in the} examples \accorpa{Freg}{F2ob} is replaced
by a more general one.
Indeed, we can \juerg{make the following assumptions:}
\Bsist
	&& F:\erre\to\errebar \quad \hbox{admits the decomposition} \quad F=\Beta +\Pi ,\quad \hbox{where}
	\qquad
	\label{hpF}
	\\[1mm]
	&& \Beta :\erre\to [0,+\infty] \quad \text{is convex and l.s.c.\
	  with subdifferential} \quad
	  \beta := \partial\Beta
	\non
	\\
	&& \quad \hbox{and fulfills} \quad 
    \Beta(0)=0
    \aand
    \lim_{r \to +\infty} {\Beta(r)}{|r|^{-2}} = +\infty \,.
    \label{hpbeta}
	\\[1mm]
	&& \Pi \in C^1(\erre) \quad \text{ with \Lip\ continuous derivative $\pi:=\Pi'$}.
	\label{hppi}
\Esist
\Accorpa\HPuno hpF hppi	
\juerg{For the source terms we assume:}
\gianni{%
\Bsist
  &&\Sph, \Ss : \erre^2 \to \erre
  \quad \hbox{are \Lip\ continuous \juerg{functions satisfying}}
  \label{hpS}
  \\
  && |\Sph(r,s)| + |\Ss(r,s)|
  \leq \Theta (|r|+|s|+1) 
  \non
  \\
  && \quad \hbox{for some constant $\Theta>0$ and every $r,s\in\erre$}\,.
  \label{hpSbis}
  \\[1mm]
  && g \in \LQ\infty\an{.}
  \label{hpg}
\Esist
\Accorpa\HPdue hpS hpg
Finally, the permeability, viscosity and sensitivity constants $\nu$, $\eta$, $\lambda$ and $\chi$ are requested to satisfy
\Bsist
  && \nu, \ \eta \in (0,+\infty)
  \aand
  \lambda, \ \chi \in [0,+\infty).
  \label{hpvisc}
\Esist
}%
\Accorpa\HPstruttura hpF hpvisc
It is \juerg{well known} that $\beta$ is a maximal 
monotone graph in $\erre\times \erre$ with corresponding domain $D(\beta)$, 
and that $0 \in \beta(0)$.
It is worth noticing that in the case $\Beta \in C^1(\erre)$ 
it follows that it is single-valued, and \oldgianni{we can write $\beta= {\Beta}'$} as well as $F'= \beta+\pi$.
Here, we immediately observe that all \juerg{of the standard} potentials \eqref{Freg}--\eqref{F2ob} fulfill \HPuno,
as well as that the biologically relevant examples given above for the source terms satisfy \last{\eqref{hpS}--\eqref{hpg}}.

\bigskip
As for the initial and boundary data, we assume:
\gianni{%
\Bsist
  && \phiz \in V
  \aand
  \Beta(\phiz) \in \Luno\,.
  \label{hpzero}
  \\	
  && \mu_\Sigma \in \H1\LdueG \cap \L2{\HxG{1/2}}\,.
  \label{hpmuSig}
  \\ 
  && \sig_\Sigma \in \LS2\,.
  \label{hpsigSig}
\Esist
}%
\Accorpa\HPdati hpzero hpsigSig
However, it is convenient to \juerg{transform} the Dirichlet inhomogeneous boundary condition for $\mu$ 
into a homogeneous one by performing a change of variable for the chemical potential.
To this end, we introduce the harmonic extension of the boundary datum $\mu_\Sigma$,
i.e., the function $h:Q\to\erre$ defined as follows\an{:}
\gianni{%
\Beq
  h(t) \in V, \quad
  \an{-}\Delta h(t) = 0
  \aand
  h(t)\suG = \an{\mu}_\Sigma(t)
  \quad \aat \,.
  \label{harm:ext}
\Eeq
We notice at once that \eqref{hpmuSig} ensures that $h$ \juerg{enjoys at least the regularity}
\begin{align}
	\gianni{h \in \H1H \cap \L2V\,.}
	\label{regh}
\end{align}
Thus, upon setting $\wt \mu:= \mu - h$,}
the boundary condition for $\mu$ in \eqref{eq:7} now becomes the homogeneous Dirichlet condition
$\wt\mu=0$ on~$\Sigma$.
However, for the sake of simplicity, we proceed with \an{abuse} of notation 
and \pier{still} denote by $\mu$ the above difference $\wt\mu$ between the chemical potential and~$h$.
This change of \an{notation} obviously \an{slightly} modifies the equations:
the \rhs\ of \eqref{eq:1} has to be adapted,
and we have to rewrite \eqref{eq:4}~as
\begin{align}
   \an{\mu 	
 \in - \badeps \,\Delta \phi
	\an{{}+ \pier{{}\badeps^{-1}} \partial F(\phi)}
  - \chi\sig
  -  h,}
  \non
\end{align}
\an{where the differential inclusion arises from the possible multi-valued 
nature of the nonlinearity $F$, \juerg{in accordance with} \eqref{hpF}--\eqref{hppi}.}
For \juerg{this reason, and in order to clarify the meaning of the equation}, 
we state the new problem \juerg{to be dealt with} 
in a precise form.
In particular, the equations \eqref{eq:1}, \eqref{eq:3} and \eqref{eq:5}\an{,} and the corresponding 
boundary conditions,
are replaced by variational equations 
(also owing to the Leibniz rule for the divergence and~\eqref{eq:2}),
and the homogeneous Dirichlet condition for $\mu$ is \juerg{enforced by} the forthcoming \eqref{regmu}.

By recalling \Stress\ and \Nutrient, and setting $m\pier{{}(\cdot)}\equiv n\pier{{}(\cdot)}\equiv 1$ 
and $\badeps=1$ as announced in the Introduction,
we look for a sextuple $\soluz$ enjoying the regularity properties
\begin{align}
	\vv & \in \L2\VV \,,
	\label{regv}
	\\
	p & \in \L{4/3}H \,,
	\label{regp}
	\\
	\mu  & \in \L2\Vz \,,
	\label{regmu}
	\\
	\phi & \in \gianni{\H1\Vzp} \cap \L\infty V \cap \L2W \,,
	\label{regphi}
	\\
	\xi & \in \L2H \,,
	\label{regxi}
	\\
	\sigma & \in \W{1,4/3}\Vp \cap \L\infty H \cap \L2V \,,
	\label{regsigma} 
\end{align}
\Accorpa\Regsoluz regv regsigma
\juerg{that satisfies} 
\begin{align}
	& \iO \gianni{\T(\vv,p)} : \nabla\bzeta 
	+ \nu \iO \vv \cdot \bzeta 
	=
	 \iO (\mu+h) \nabla\phi \cdot \bzeta
	+ \iO \Ns(\phi,\sig) \nabla\sig \cdot \bzeta
	\non
	\\
	& \quad \hbox{for every $\bzeta\in\VV$ and \aet}\,,
	\label{prima}
	\\[2mm]
	\separa
	& \div\vv  
	= g 
	\quad \aeQ\,,
	\label{eqdiv}
	\\[1mm]
	\separa
	& \< \dt \phi, \badphi>_{\Vz}
	+ \iO \nabla\mu \cdot \nabla\badphi
	= \iO \Sph(\phi,\sig) \badphi
	- \iO (\nabla\phi \cdot \vv
	+ \phi g\an{)} \, \badphi
	\non
	\\
	& \quad \hbox{for every $\badphi\in\Vz$ and \aet}\,,
	\label{seconda}
	\\[1mm]
	\separa
	& \iO \nabla\phi \cdot \nabla z
	+ \iO (\xi + \pi(\phi)) z
	= \iO \bigl( \mu + h - \Nph(\phi,\sig) \bigr) z
	\non
	\\
	& \quad \hbox{for every $z\in V$ and \aet}\,,
	\label{terza}
	\\[2mm]
	\separa
	& \xi \in \beta(\phi)
	\quad \aeQ \,,
	\label{xibetaphi}
	\\[1mm]
	\separa
	& \< \dt \sig, \zeta>_V
	+ \iO \nabla\Ns(\phi,\sig) \cdot \nabla\zeta
	\non
	\\
	& = \iO \Ss (\phi,\sig) \zeta
	- \iO (\nabla\sig \cdot \vv + \sig g)\zeta
	+ \kappa \iG (\sig_\Sigma - \sig) \zeta 
	\non
	\\
	& \quad \hbox{for every $\zeta\in V$ and \aet}\,,
	\label{quarta}
\end{align}
as well as the initial conditions
\Beq
	\pcol{\phi(0) = \phiz
	\aand
	\sig(0) = \sigz \,.}
	\label{cauchy}
\Eeq
\Accorpa\Pbl prima cauchy

Here is our \an{main} result:

\Bthm
\label{Existence}
Assume \HPstruttura, and let the notations \Stress\ and \Nutrient\ be in force.
Moreover, let \HPdati\ be fulfilled,
and let $h$ be defined by \eqref{harm:ext}.
Then, the weak formulation \Pbl\ of the Cahn--Hilliard--Brinkman system admits at least one solution 
\an{$\soluz$} with the regularity specified by \Regsoluz.
\Ethm

\Brem
\label{Morereg}
\an{It is worth pointing out that equation \eqref{terza} has been formulated in a} weak form \an{just} for convenience.
\an{Indeed}, thanks to \an{\eqref{xibetaphi}, along with the regularity properties \eqref{regphi} and \eqref{regxi},}
the variational equation \eqref{terza} is equivalent to the boundary value problem
\Beq
  - \Delta\phi + \xi + \pi(\phi)
  = \mu + h - \Nph(\phi,\sig)
  \enskip \aeQ 
  \aand
  \dn\phi=0
  \enskip \hbox{on $\Sigma$} \,.
  \label{strongterza}
\Eeq
As \an{further} regularity is concerned,
the components $\phi$ and $\xi$ of every solution satisfy
\Beq
  \phi \in \L2{\Wx{2,6}} \cap \L4\Hdue
  \aand
  \xi \in \L2{\Lx6}
  \label{morereg}
\Eeq
(in~the two-dimensional case the summability exponent $6$ can be replaced by any $q\geq1$),
as shown in the forthcoming Remark~\ref{Regphi}.
Finally, we notice that, \juerg{in view of the very low regularity at disposal for the nutrient variable $\sig$ and \an{the} velocity field~$\vv$, the uniqueness of weak solutions is not to be expected}.
\Erem

\gianni{We continue this section by listing some tools \an{that will be useful later on}.
We first recall Young's inequality
\Beq
  \an{a\,b \leq \frac \delta q \, a^{q} + \frac{{(\delta)}^{-q'/q}}{q'} \, b^{q'}}
  \quad \hbox{for all $a,b\in[0,+\infty)$, $q\in(1,+\infty)$ and $\delta>0$},
  \label{young}
\Eeq
\juerg{where $q'$ denotes the conjugate exponent of $q$ given by the identity} $\,(1/q)+(1/q')=1$.
We repeatedly use it, mainly with $q=\an{{}q'={}}2$.
We also account for H\"older's inequality,
as well as for the following Sobolev, compactness, Poincar\'e, and Korn inequalities:
\begin{alignat}{2}
  & \norma v_q
  \leq \gianni{C_S} \, \normaV v
  \quad && \hbox{for every $v\in V$ and $q\in[1,6]$}\,,
  \label{sobolev}
  \\[2mm]
  & \norma v_4^2
  \leq \delta \, \norma{\nabla v}^2 + C_\delta \, \norma v^2
  \quad && \hbox{for every $v\in V$ and $\delta>0$}\,,
  \label{compact}
  \\[2mm]
  & \normaV v
  \leq \gianni{C_P} \, \norma{\nabla v}
  \quad && \hbox{for every $v\in\Vz$}\,,
  \label{poincare}
  \\[2mm]
  & \normaV\vv^2
  \leq \gianni{C_K} \bigl(
    \norma\vv^2
    + \norma{D\vv}^2
  \bigr)
  \quad && \an{\hbox{for every $\vv\in\VV$}}\,.
  \label{korn}
\end{alignat}
\juerg{Here,} \pier{the constants} $C_S$, $C_P$, and $C_K$, depend only on $\Omega$,
while $C_\delta$ depends on~$\delta$, in addition.
In \eqref{korn}, the notation~\eqref{symgrad} is used}.

\gianni{Next, we present three \an{auxiliary} results \juerg{that will be used} in the sequel.
The first one is stated in a more general setting 
as an exercise in \cite[Ex.~III.3.5]{Galdi},
but it \an{readily follows as a} corollary \juerg{from} \cite[Thm.~III.3.1]{Galdi}.
The second one \an{is related to the Stokes resolvent operator and} is a particular case of \cite[\an{Thm.~3}]{AbelsTer} \pier{(which is an extension of the results in \cite{Giga}). Finally,} the last one regards the trace operator.}

\begin{lemma}
\label{Lemmau}
There exists a constant~$C$ that depends only on $\Omega$ such that
for every $f$ and $\aa$ satisfying 
\Beq
  f \in \an{H} , \quad
  \aa \in \HHxG{1/2},
  \aand
  \iO f = \iG \aa\cdot\nn\,,
  \label{compatib}
\Eeq
there exists \juerg{some} $\uu\in\an{\VV}$ satisfying
\Beq
  \div\uu = f \quad \aeO, \quad
  \uu\suG=\aa,
  \aand
  \norma\uu_{\an{V}}
  \leq C \bigl(
    \norma f + \norma\aa_{\HHxG{1/2}}
  \bigr) \,.
  \label{diveq}
\Eeq
\end{lemma}

\Blem
\label{AbTe}
Assume that $\ff\in\an{\HH}$ and $f\in\an{V}$.
Then, there exists a unique pair $(\vv,p)$ satisfying
\Bsist
  && \vv \in \HHx2 
  \aand
  p \in \an{V}\,,
  \non
  \\[2mm]
  && - \div\T(\vv,p) + \nu\vv = \ff 
  \aand
  \div\vv = f
  \quad \hbox{in $\Omega$}\,,
  \non
  \\[2mm]
  && \an{\T(\vv,p)}\nn = \0 
  \quad \hbox{on $\Gamma$}\,.
  \non
\Esist
Moreover\an{,} the \juerg{mapping} $\,\Psi=(\Psi_1,\Psi_2):(\ff,f)\mapsto(\vv,p)$\, is linear and continuous 
from $\an{\HH}\times\an{V}$ into $\HHx2\times\an{V}$.
\Elem

\Blem
\label{Traccia}
The trace operator maps $\L\infty H\cap\L2V$ into $\L4\LdueG$,
and \juerg{it holds the estimate}
\Beq
  \ioT \norma{v(t)}_{\LdueG}^4 \, dt
  \leq C \, \norma v_{\L\infty H}^2 \, \norma v_{\L2V}^2 
  \label{traccia}
\Eeq
for every $v\in\L\infty H\cap\L2V$,
where the constant $C$ depends only on $\Omega$,
and \juerg{where} $\,v\,$ also denotes the trace of~$\,v\,$ on $\an{\Gamma}$.
\Elem

\begin{proof}[\last{Proof of Lemma \ref{Traccia}}]
As we did not find a precise reference,
we provide a sketch of the proof.
We denote by $C_1$, $C_2$,\dots\ constants that depend only on~$\Omega$.
We introduce the \an{real} interpolation space
(by the way, the Besov space $B^{1/2}_{2,1}(\Omega)$)
\Beq
  B := (V,H)_{1/2,1} \,.
  \label{besov}
\Eeq
The trace operator $v\mapsto v\suG$ maps $B$ into $\LdueG$
and is linear and continuous.
In the \juerg{half-space} case $\Omega=\last{\erre^3_+}$, 
this can be deduced, e.g., from formula (I.17) 
(with the notation $Y(1,\last{\erre^3_+})$ for~\eqref{besov})
in~the paper \cite{Tartar},
where it is also shown that the operator maps $B$ onto $\LdueG$.
\juerg{As usual, the result is then} extended to the general case by \juerg{using local charts and a} partition of unity.
This leads to the estimate
\Beq
  \norma v_{\LdueG} \leq C_1 \, \norma v_B
  \quad \hbox{for every $v\in B$}.
  \non
\Eeq
On the other hand, the interpolation inequality
\Beq
  \norma v_B \leq C_2 \, \normaV v^{1/2} \, \norma v^{1/2}
  \non
\Eeq
holds true for every $v\in V$.
Now, \pier{letting} $v\in\L\infty H\cap\L2V$\an{, we} have \aat\ that
\Bsist
  && \norma{v(t)}_{\LdueG} 
  \,\leq C_1\, \norma{v(t)}_B \,\leq C_3 \, \normaV{v(t)}^{1/2} \, \norma{v(t)}^{1/2}
	\,\leq\, C_3 \, \normaV{v(t)}^{1/2} \, \norma v_{\L\infty H}^{1/2}\,,
  \non
\Esist
and \eqref{traccia} \an{directly} follows by taking the \juerg{$4{th}$} powers and integrating over $(0,T)$.
\end{proof}

\Brem
\an{We notice that the application of the trace estimate in \cite[Thm.~II.4.1]{Galdi} (whose proof is left to the reader), with the parameters therein being chosen
as $r = q = 2, m = 1, n = d=3, \lambda = 0$, produces a similar trace inequality\juerg{, namely,}
\begin{align*}
	\norma{v}_{\LdueG} \leq c\,( \norma{v} + \norma{v}^{1/2}\norma{v}^{1/2}_{V})
	\quad 
	\text{for every $v \in V.$}
\end{align*}}%
\Erem

\an{Besides, let us state} a general rule concerning the constants 
that appear in the estimates to be performed in the following.
The small-case symbol $\,c\,$ stands for a generic constant
whose actual \an{value} may change from line to line, and even within the same line,
and depends only on~$\Omega$, \an{the} shape of the nonlinearities,
and \an{the} constants and the norms of the functions involved in the assumptions of the statements.
In particular, the values of $\,c\,$ do not depend on the parameters $\,\eps\,$ and $\,k\,$ \juerg{that will
be introduced} in the next section.
A~small-case symbol with a subscript like $c_\delta$ (\an{specifically}, with $\delta=\eps$)
indicates that the constant may depend on the parameter~$\delta$, in addition.
On the contrary, we mark precise constants that we can refer~to
by using different symbols
(see, e.g., \eqref{sobolev} and \eqref{diveq}).

\an{%
The next sections aim to rigorously prove Theorem~\ref{Existence}.
A standard approach to guarantee the existence of solutions to similar Cahn--Hilliard type systems is based on suitable approximation procedures that can be schematized as follows: first\pier{,} one has to regularize the possibly singular nonlinearity $F$, with the classical choice being the Yosida regularization that \pier{depends} on a parameter, say, 
$\eps>0$. Then, for every \juerg{fixed $\eps>0$}, one further discretizes the system in space 
via the Galerkin method and solves
a family of finite-dimensional problems that depend on a parameter \juerg{$k\in\enne$}.
 Next, one has to provide some rigorous estimates, \juerg{which are} independent of both $\,\eps\,$ and $\,k$.
From \juerg{these estimates}, using  weak and weak star compactness arguments, one can find a suitable subsequence and eventually pass to the limit as $k \to \infty$ and as $\eps \to 0$\pier{, thus showing} that the obtained limits yield a solution to the original system.}

\an{%
The essence of the proof of Theorem~\ref{Existence} can be roughly schematized by the abovementioned steps, but we had to face a major obstacle arising from the different boundary conditions in \eqref{eq:7}.
\juerg{The} latter prevent the solvability of the ODE system \juerg{originating} from the Galerkin scheme. 
In fact, \juerg{due to \eqref{eq:7}, one is in the approximation naturally led} to consider Schauder bases (cf. \eqref{eigenj} and \eqref{eigenzj}) for the Laplace operator 
with homogeneous Neumann and Dirichlet boundary conditions for \pier{$\phi, \ \sigma$ and} 
$\mu$, respectively. 
Despite \juerg{of} being natural, this choice completely impedes the solvability of the discrete problem, \juerg{because there
occur inner products between the elements of the different two \pier{bases}}.
To overcome this intrinsic \juerg{difficulty}, we introduce an intermediate approximation step, which consists in adding 
\juerg{further regularizing terms} at the level $\eps>0$ \juerg{that somehow dominate the mixed terms and enable us}
to solve the Galerkin system (cf. \eqref{primak}--\eqref{cauchyk}).
Finally, we pass to the limit as $\eps \searrow 0$ as anticipated above\juerg{, thus proving} the theorem.}%


\section{Approximation}
\label{APPROX}
\setcounter{equation}{0}

In this section, we introduce and solve a proper approximating problem depending on the parameter~$\eps>0$.
First of all, we replace the functional $\Beta$ and the graph $\beta$ by their Moreau\an{--}Yosida regularizations
$\Betaeps$ and~$\betaeps$, respectively
(see, e.g., \cite[pp.~28 and~39]{Brezis}).
Then, \an{we set}
\Beq
  \Feps := \Betaeps + \Pi
  \label{defFeps}
\Eeq
\an{and recall that classical theory of convex analysis entails \juerg{the following facts}:}
\Bsist
  && \hbox{$\betaeps$ is monotone and \Lip\ continuous with $\betaeps(0)=0$}\,,
  \label{monbetaeps}
  \\
  && 0 \leq \Betaeps(r) = \int_0^r \betaeps(s) \, ds \leq \Beta(r)
  \quad \hbox{for every $r\in\erre$}\,,
  \label{disugBetaeps}
  \\
  && \hbox{for every $M>0$ there exist $C_M>0$ and $\eps_M>0$ such that}
  \non
  \\
  && \quad 
  \Feps(r) 
  \geq M \, r^2 - C_M
  \quad \hbox{for every $r\in\erre$ and every $\eps\in(0,\eps_M)$}\,.
  \qquad
  \label{coerceps}
\Esist
\an{The only nonobvious fact is the coercivity property in \eqref{coerceps}. In this direction, let us 
fix} $M>0$ and observe that our assumption \eqref{hpbeta} on $\Beta$ implies that
\Beq
  \an{\Beta}(r)
  \geq 2M \, r^2 - C_M \quad \hbox{for every $r\in\erre$ and some constant $C_M>0$}.
  \non
\Eeq
It \an{then} follows that
\Bsist
  && \Betaeps(r)
  := \inf_{s\in\erre} \left\{ \frac 1 {2\eps} |s-r|^2 + \Beta(s) \right\} 
  \non
  \\
  && \geq \inf_{s\in\erre} \left\{ \frac 1 {2\eps} |s-r|^2 + 2M \, s^2 - C_M \right\} 
  = \frac 1 {2\eps} |s_*-r|^2 + 2M \, s_*^2 - C_M\,,
  \non
\Esist
where $s_*$ is the minimum point, namely, $s_*=r/(1+4M\eps)$.
Hence, we have that
\Beq
  \Betaeps(r)
  \geq 2M \, s_*^2 - C_M
  = \frac {2M} {(1+4M\eps)^2} \, r^2 - C_M
  \geq M \, r^2 - C_M
  \non
\Eeq
for every $r\in\erre$, whenever $(1+4M\eps)^2\leq2$,
i.e., \eqref{coerceps} with $\Feps$ replaced by~$\Betaeps$ (with an obvious choice of~$\eps_M$).
Then, \eqref{coerceps} itself follows, since $\pi$ is \Lip\ continuous.

Besides this regularization, we replace $g$ and $h$ by smoother functions $\geps$ and $\heps$ satisfying
\Bsist
  && \geps \in \CQ0 , \quad
  \norma\geps_\infty \leq c \,,
	\aand
  \geps \to g
  \quad \hbox{\aeQ\ as $\eps\searrow 0$}\,,
  \label{hpgeps}
  \\[2mm]
  && \heps \in \H1H \cap \L2V \cap \CQ0 \,, \quad
  \norma\heps_{\H1H\cap\L2V} \leq c \,,\non\\
  &&\quad \hbox{\an{and $\heps\to h$ strongly} in $\L2V$ as $\eps\searrow 0$\pier{.}}
  \label{hpheps}
\Esist
\an{For simplicity, we do not enter in the details concerning the construction of the regularizations above. 
However, let us mention that standard mollification arguments are enough to get the desired properties prescribed by \eqref{hpgeps} and \eqref{hpheps}, respectively.}
\pier{Moreover, as anticipated, we introduce artificial viscous terms in some of the equations.} We prefer to present all of them in their variational form.
The approximating problem thus consists in finding a quintuple $\soluzeps$
that satisfies the regularity properties
\begin{align}
	\vveps & \in \L2\VV \,,
	\label{regveps}
	\\
	\peps & \in \L{4/3}H \,,
	\label{regpeps}
	\\
	\mueps  & \in \H1\Vzp \last{{}\cap \L\infty H{}}\cap \L2\Vz \,,
	\label{regmueps}
	\\
	\phieps & \in \H1H \cap \L\infty V \cap \L2W \,,
	\label{regphieps}
	\\
	\sigeps & \in \W{1,4/3}\Vp \cap \L\infty H \cap \L2V\,,
	\label{regsigeps} 
\end{align}
\Accorpa\Regsoluzeps regveps regsigeps
and solves the variational equations
\begin{align}
	& \iO \T(\vveps,\peps) : \nabla\bzeta 
	+ \nu \iO \vveps \cdot \bzeta 
	=
	 \iO (\mueps+\heps) \nabla\phieps \cdot \bzeta
	+ \iO \Ns(\phieps,\sigeps) \nabla\sigeps \cdot \bzeta
	\non
	\\
	& \quad \hbox{for every $\bzeta\in\VV$ and \aet}\,,
	\label{primaeps}
	\\[2mm]
	\separa
	& \div\vveps
	= \geps
	\quad \aeQ
	\label{eqdiveps}
	\\[1mm]
	& \< \dt(\eps\mueps+\phieps) , \badphi >_{\Vz}
	+ \iO \nabla\mueps \cdot \nabla\badphi
	= \iO \Sph(\phieps,\sigeps) \badphi
	- \iO (\nabla\phieps \cdot \vveps + \phieps\geps) \badphi
	\non
	\\
	& \quad \hbox{for every $\badphi\in\Vz$ and \aet}\,,
	\label{secondaeps}
	\\
	\separa
	& \eps \iO \dt\phieps \, z 
    + \iO \nabla\phieps \cdot \nabla z
	+ \iO \Feps'(\phieps) \, z
	= \iO \bigl( \mueps \an{{}+ \heps{}} - \Nph(\phieps,\sigeps) \bigr) z
	\non
	\\
	& \quad \hbox{for every $z\in V$ and \aet}\,,
	\label{terzaeps}
	\\
	\separa
	& \< \dt \sigeps, \zeta>_V
	+ \iO \nabla\Ns(\phieps,\sigeps) \cdot \nabla\zeta
	\non
	\\
	& = \iO \Ss (\phieps,\sigeps) \zeta
	- \iO (\nabla\sigeps \cdot \vveps + \sigeps \geps) \zeta
	+ \kappa \iG (\sig_\Sigma - \sigeps) \zeta 
	\non
	\\
	& \quad \hbox{for every $\zeta\in V$ and \aet}\,,
	\label{quartaeps}
\end{align}
as well as the initial conditions
\Beq
	\pcol{\phieps(0) = \phiz \,, \quad
    \mueps(0) = 0
	\aand
	\sigeps(0) = \sigz \,.}
	\label{cauchyeps}
\Eeq
\Accorpa\Pbleps primaeps cauchyeps

\an{Let us incidentally notice that the presence of a selection $\xi$ (cf. \eqref{xibetaphi}) is no longer needed in the approximating $\eps$-problem due to the regularity at disposal for $F_\eps$ defined by \eqref{defFeps}.}
\pier{Moreover, we can simply write $F'_\eps(\phi)$ in place of  $\betaeps(\phi) +\pi (\phi) $ in \eqref{terzaeps}.}
\Brem
\label{Regmueps}
We notice that, due to \eqref{regphieps} and the initial condition for~$\phieps$, 
the regularity for $\dt\mueps$ given in \eqref{regmueps} and the initial conditions for $\mueps$
are equivalent~to
\Beq
  \dt(\eps\mueps+\phieps) \in \L2\Vzp
  \aand
  (\eps\mueps+\phieps)(0) = \phiz\,,
  \non
\Eeq
respectively.
\Erem

\Bthm
\label{Existenceeps}
Under the assumptions of Theorem~\ref{Existence} and with the above notation,
the approximating problem \an{{}\Pbleps{}} has at least one solution \an{$\soluzeps$ which fulfills \eqref{regveps}--\eqref{regsigeps}}.
\Ethm

The rest of \juerg{this} section is devoted to the proof of this theorem.
The method we use starts from a discrete problem based on a Faedo\an{--}Galerkin scheme.
To this end, we introduce the nondecreasing sequences 
$\{\lambdaj\}$ and $\{\lambdazj\}$ of \an{ eigenvalues}
and the corresponding complete orthonormal sequences $\{\ej\}$ and $\{\ezj\}$ of \an{ eigenfunctions}
of the eigenvalue problems for the Laplace operator 
with homogeneous Neumann and Dirichlet boundary conditions, respectively.
Namely, we have~that
\an{
\begin{alignat}{2}
  & - \Delta\ej = \lambdaj\ej
  \quad \hbox{in $\Omega$}
  \aand
   \dn\ej = 0
  \quad && \hbox{on $\Gamma$}\,,
  \label{eigenj}
  \\
  & - \Delta\ezj = \juerg{\lambdazj}\ezj
  \quad \hbox{in $\Omega$}
  \aand
  \ezj = 0
  \quad && \hbox{on $\Gamma$}\,,
  \label{eigenzj}
\end{alignat}}
for $j=1,2,\dots$, as well as \an{the normalization conditions}
\Beq
  \iO \ei \ej 
  = \iO \ezi \ezj
  = \delta_{ij} 
  \quad \hbox{for every $i$ and $j$} \,\an{,}
  \label{orthonormal}
\Eeq
\juerg{with the standard Kronecker symbol $\delta_{ij}$}.
Moreover, if we set, for $k=1,2,\dots,$
\Beq
  \Vk := \Span\{\ej:\ 1\leq j\leq k\}
  \aand
  \Vzk := \Span\{\ezj:\ 1\leq j\leq k\}	\an{\,,}
  \label{defspazik}
\Eeq
\juerg{then} the unions of the these spaces are dense in $V$ and~$\Vz$, respectively,
and both are dense in $H$ as well.
We notice at once that all \juerg{of} the above eigenfunctions are smooth since $\Omega$ is smooth\an{,}
that $\,\lambda_1=0$\an{, and that $\,e_1=|\Omega|^{-1/2}$}.
Then, the discrete problem related to $k$ consists in finding a \an{quintuple}
$\soluzk$ with the regularity specified~by
\Bsist
  && \vvk \in \L2{\HHx2} , \quad
  \pk \in \L2{\an{V}} , \quad
  \phik, \sigk \in C^1([0,T);\Vk) \cap \LQ\infty,\qquad
  \non
  \\
  && 
  \aand
  \muk \in C^1([0,T);\Vzk) \cap \LQ\infty, 
  \qquad
  \label{regsoluzk}
\Esist
that solves the system
\begin{align}
	& \iO \T(\vvk,\pk) : \nabla\bzeta 
	+ \nu \iO \vvk \cdot \bzeta 
	=
	 \iO (\muk + \heps) \nabla\phik \cdot \bzeta
	+ \iO \Ns(\phik,\sigk) \nabla\sigk \cdot \bzeta\,,
	\label{primak}
	\\[2mm]
	& \div\vvk
	= \geps
	\quad \aeQ \,,
	\label{eqdivk}
	\\[2mm]
	& \eps \iO \dt\muk \, \badphi
	+ \iO \nabla\muk \cdot \nabla\badphi
	\non
	\\
	& = - \iO \dt\phik \, \badphi
	+ \iO \Sph(\phik,\sigk) \badphi
	- \iO (\nabla\phik \cdot \vvk + \phik\geps) \badphi\,,
	\label{secondak}
	\\[1mm]
	& \eps \iO \dt\phik \, z 
    + \iO \nabla\phik \cdot \nabla z
	+ \iO \Feps'(\phik) \, z
	= \iO \bigl( \muk \an{{}+ \heps{}} - \Nph(\phik,\sigk)\bigr) z\,,
	\label{terzak}
	\\[1mm]
	& \iO \dt \sigk \, \zeta
	+ \iO \nabla\Ns(\phik,\sigk) \cdot \nabla\zeta
	\non
	\\
	& = \iO \Ss (\phik,\sigk) \zeta
	- \iO (\nabla\sigk \cdot \vvk + \sigk \geps) \zeta
	+ \kappa \iG (\sig_\Sigma - \sigk) \zeta \,,
	\label{quartak}
\end{align}
for every $\bzeta\in\VV$, $\badphi\in\Vzk$, $z,\,\zeta\in\Vk$, and $t\in[0,T)$,
and that fulfills the initial conditions
\Beq
	\iO \phik(0) \badphi = \iO \phiz \badphi \,, \enskip
	\muk(0) = 0\, ,
	\aand
	\iO \sigk(0) \zeta = \iO \sigz \zeta
	\label{cauchyk}
\Eeq
for every $\badphi\in\Vzk$ and $\zeta\in\Vk$.
\Accorpa\Pblk primak cauchyk

\step
Existence \pier{for} the discrete problem

The first aim of ours is to show the existence of at least one solution
(we~do not care about uniqueness, since it is not needed).
\pier{The} method relies on a proper application of Lemma~\ref{AbTe}.
\an{Besides, let us point out that the idea of employing the Stokes resolvent to express the velocity field $\vv_k$ in terms of the other variables $\phi_k,\mu_k,$ and $\sigma_k$ is largely inspired by \cite{EGAR} (see also \cite{EGAR2, KS2}).}
For a while, the symbols $\phik$, $\muk$ and $\sigk$ denote independent variables. 
To every triplet $(\phik,\muk,\sigk)\in\Vk\times\Vzk\times\Vk$\an{,}
we associate the \an{vector-valued} function 
\Beq
  \ff_k 
  := (\muk+\heps) \nabla\phik
  + (\sigk+\chi(1-\phik)) \nabla\sigk\,,
  \label{ffk}
\Eeq
and we notice that $\ff_k$ only depends on time through the continuous function~$\heps$\juerg{,
while the dependence} on space occurs just through the eigenfunctions and their gradients, 
which are smooth.
In particular, on \juerg{the one hand}, 
if we read $\ff_k$ as a function of the coefficients of $\phik$, $\muk$\an{,} $\sigk$ 
(with respect to the bases just chosen)\an{,} and~$t$,
\juerg{then} we see that it is continuous.
On the other \juerg{hand}, for every $t\in[0,T]$, we are allowed to apply Lemma~\ref{AbTe} with
$\ff=\ff_k(t)$ and $f=\geps(t)$.
Since the \juerg{mapping} $\Psi$ is linear, continuous, and time independent, and \juerg{since} $\geps$ and $\heps$ are continuous,
this yields a pair of functions 
\an{that} are continuous with respect to the coefficients of $\phik$, $\muk$, $\sigk$, and~$t$.
This observation is made to ensure the continuity of the functions that rule the system of ODE's we are going to introduce.

Now, we let $\phik$, $\muk$ and $\sigk$ depend on time.
To every triplet $(\phik,\muk,\sigk)\in\L\infty{\Vk\times\Vzk\times\Vk}$ 
we associate the function $\ff_k$ still given by \eqref{ffk} and, for every $t\in[0,T]$, 
we apply Lemma~\ref{AbTe} as before.
We obtain two functions, which we still term
$\Psi_1(\ff_k,\geps)$ and $\Psi_2(\ff_k,\geps)$ with an abuse of notation,
that\an{, according to the lemma,} belong to $\L\infty{\HHx2}$ and $\L\infty{\last{{}V{}}}$, respectively.
By construction, the pair of functions 
$(\vvk,\pk):=(\Psi_1(\ff_k,\geps),\Psi_2(\ff_k,\geps))$ solves the equations \accorpa{primak}{eqdivk} 
corresponding to the given triplet $(\phik,\muk,\sigk)$.
Therefore, the whole problem \Pblk\ is equivalent to the problem of finding
a triplet $(\phik,\muk,\sigk)$ with the regularity specified in \eqref{regsoluzk}
that solves
\begin{align}
	& \eps \< \dt \muk, \badphi>_{\an{V_0}}
	+ \iO \nabla\muk \cdot \nabla\badphi
	\non
	\\
	& = - \< \dt \phik, \badphi>_{\an{V_0}}
	+ \iO \Sph(\phik,\sigk) \badphi
	- \iO (\nabla\phik \cdot \Psi_1(\ff_k,\geps) + \phik\geps) \badphi\,,
	\label{secondakbis}
	\\
	& \eps \< \dt\phik , z >_V
    + \iO \nabla\phik \cdot \nabla z
	+ \iO \Feps'(\phik) \, z
	= \iO \bigl( \muk \an{{}+ \heps{}} - \Nph(\phik,\sigk)  \bigr) z\,,
	\label{terzakbis}
	\\
	& \< \dt \sigk, \zeta>_V
	+ \iO \nabla\Ns(\phik,\sigk) \cdot \nabla\zeta
	\non
	\\
	& = \iO \Ss (\phik,\sigk) \zeta
	- \iO (\nabla\sigk \cdot \Psi_1(\ff_k,\geps) + \sigk \geps) \zeta
	+ \kappa \iG (\sig_\Sigma - \sigk) \zeta \,,
	\label{quartakbis}
\end{align}
for every $\badphi\in\Vzk$, $z,\,\zeta\in\Vk$ and $t\in[0,T)$,
with $\ff_k$ given by \eqref{ffk},
and satisfies the initial conditions \last{in} \eqref{cauchyk}.
We show that this problem has at least one solution.
To~this end, we represent the unknowns in terms of the bases of the spaces $\Vk$ and $\Vzk$,~i.e.,
\Beq
  \phik(t) = \somma j1k \phikj(t) \ej , \quad
  \muk(t) = \somma j1k \an{\mukj}(t) \ezj \,,
  \aand
  \sigk(t) = \somma j1k \an{\sigkj}(t) \ej \,,
  \non
\Eeq
and introduce the $\erre^k$-valued functions
\Beq
  \hphik := (\phikj)_{j=1}^k \,, \quad
  \hmuk := (\mukj)_{j=1}^k \,,
  \aand
  \hsigk := (\sigkj)_{j=1}^k \,,
  \non
\Eeq
which are the true unknowns.
In terms of these coefficients, the discrete problem takes the form
\Bsist
  && \an{\eps} \hmuk'(t) 
  = \an{\AA}_{\eps,k}(\hphik(t),\hmuk(t),\hsigk(t),t)
  - \hphik'(t)\,,
  \label{secondakter}
  \\
  && \an{\eps} \hphik'(t) 
  = \an{\BB}_{\eps,k}(\hphik(t),\hmuk(t),\hsigk(t),t)\,,
  \label{terzakter}
  \\
  && \hsigk'(t) 
  = \an{\CC}_{\eps,k}(\hphik(t),\hmuk(t),\hsigk(t),t)\,,
  \label{quartakter}
\Esist
with some continuous functions 
$\an{\AA}_{\eps,k},\an{\BB}_{\eps,k},\an{\CC}_{\eps,k}:\erre^k\times\erre^k\times\erre^k\times[0,T]\to\erre^k$,
and the initial conditions for $\hphik$, $\hmuk$ and $\hsigk$ are trivially derived from~\eqref{cauchyk}.
By replacing $\hphik'$ in \eqref{secondakter} \juerg{using~\eqref{terzakter} (recall that now $\eps>0$ is fixed),} 
we obtain a standard Cauchy problem
for a $3k$-dimensional nonlinear \juerg{ODE} system ruled by a continuous function.
This allows us to apply the \an{Cauchy--}Peano theorem,
which ensures the existence of at least one local solution.
This local solution can be extended to a maximal solution,
which provides a maximal solution $(\phik,\muk,\sigk)$ to \Pblk\
defined in the interval $[0,T_k)$ for some $T_k\in(0,T]$.
We claim that this solution is bounded (as required) and global, i.e., that $T_k=T$.
The proof relies on the estimate
\Beq
  \norma\phik_{L^\infty(0,T_k;H)}
  + \norma\muk_{L^\infty(0,T_k;H)}
  + \norma\sigk_{L^\infty(0,T_k;H)}
  \leq c_\eps\,,
  \label{fromprimastima}
\Eeq
\juerg{which we are going to} prove in the next lines.
Since \eqref{orthonormal} implies that
\Beq
  \norma{\phik(t)}^2
  = \somma j1k |\phikj(t)|^2
  = |\hphik(t)|^2
  \quad \hbox{for every $t\in[0,T_k)$}\,,
  \non
\Eeq
and similarly for the other two components,
\eqref{fromprimastima} shows that the $\erre^{3k}$-valued function $(\hphik,\hmuk,\hsigk)$ is bounded.
Then, maximality also implies that the solution is global.

\medskip

\juerg{Next}, we perform a number of a priori estimates that \juerg{will} allow us to let $k$ tend to infinity
and to show that the approximating problem \Pbleps\ has at least one solution.
The first of \juerg{these estimates proves the validity of}~\eqref{fromprimastima}, in particular.

\Brem
\label{Projk}
We note that the initial values $\phik(0)$ and $\sigk(0)$ 
are the $H$-projections of $\phiz$ and $\sigz$ onto~$\Vk$.
To deal with them and for other purposes, it is worth noting a property 
of the $H$-projection $\vk$ of a generic element $v\in V$ onto~$\Vk$.
We have that $\norma\vk\leq\norma v$, and we derive a similar inequality for the gradients.
By definition, for $j=1,\dots,k$ we have~that  
\Beq
  \iO \vk \ej = \iO v \ej\,,
  \non
\Eeq
whence\an{, using \eqref{eigenj},} also
\Bsist
  && \iO \nabla\vk \cdot \nabla\ej
  = - \iO \vk \Delta\ej
  = - \lambdaj \iO \vk \ej
  \non
  \\
  && = - \lambdaj \iO v \ej
  = - \iO v \Delta\ej
  = \iO \nabla v \cdot \nabla\ej \,.
  \non
\Esist
By linear combination, it follows that
\Beq
  \iO \nabla\vk \cdot \nabla w
  = \iO \nabla v \cdot \nabla w
  \quad \hbox{for every $w\in\Vk$}\,,
  \non
\Eeq
and we conclude that 
\Beq
  \hbox{the $H$-projection $\vk$ coincides with the $V$-projection of~$v$}.
  \non
\Eeq
By noting that the choice $w=\vk$ is admitted in the above identity,
and collecting everything, we conclude that
\Beq
  \norma\vk \leq \norma v , \quad
  \normaV\vk \leq \normaV v,
  \aand
  \norma{\nabla\vk} \leq \norma{\nabla v},
  \quad \hbox{for every $v\in V$}.
  \label{genprojk}
\Eeq
Therefore, in particular, we \juerg{have} the inequalities
\Beq
  \norma{\phik(0\an{)}}\leq\norma\phiz\,, \quad
  \norma{\sigk(0\an{)}}\leq\norma\sigz\,,
  \aand
  \norma{\nabla\phik(0\an{)}}\leq\norma{\nabla\phiz} \,.
  \label{projk}
\Eeq
\juerg{Now let} $q\in[1,+\infty]$ and $v\in\L qV$, and define $\vk:Q\to\erre$ \juerg{as follows}:
\aat, $\vk(t)$ is the $H$-projection of $v(t)$ onto~$\Vk$.
Then, 
\Beq
  \vk \in \L q\Vk\,,
  \aand
  \norma\vk_{\L qV} \leq \norma v_{\L qV} \,.
  \non
\Eeq
Moreover, if $q<+\infty$, \juerg{then} we also have that
$\vk\to v$ strongly in $\L qV$.
Indeed,
\Beq
  \vk(t) \to v(t)
  \quad \hbox{strongly in $V$}\,,
  \aand
  \normaV{\vk(t)-v(t)}^q \leq 2^q \, \norma{v(t)}^q
  \quad \aat \,,
  \non
\Eeq
so that one can apply the Lebesgue dominated convergence theorem.
Clearly, everything can be repeated for the space $\Vz$ and the projection on~$\Vzk$.
\Erem

\step
First a priori estimate

We recall that we do not \juerg{yet know that $T_k=T$}.
For every $t\in(0,T_k)$, we apply Lemma~\ref{Lemmau} \an{with} the choices
\Beq
  \an{f} = \geps(t)
  \aand
  \aa = \frac 1{|\Gamma|} \, \Bigl( \textstyle\iO \geps(t) \Bigr) \nn\,,
  \non
\Eeq
by observing that the assumptions \eqref{compatib} are satisfied,
and we term $\uuk(t)$ the function given by the lemma.
\an{Avoiding writing} the time $t$ for a while\juerg{, and} recalling \eqref{diveq} and \eqref{hpgeps}, we have that
\Beq
  \norma\uuk_{\an{L^\infty(0,T_k; \pier{{}\VV{}})}} \leq c \,.
  \label{disuguk}
\Eeq
Then, we test \eqref{primak} by $\bzeta=\vvk-\uuk$, \an{with the aim of removing the pressure from the first estimate.
This idea has been introduced in \cite{EGAR}, and it relies on the identities $D\vvk:\nabla\vvk=|D\vvk|^2$ and $\I:\nabla(\vvk-\uuk)=\div(\vvk-\uuk)=0$.
Besides,} by also integrating over $(0,t)$ with respect to time, 
where $t\in(0,T_k)$ is arbitrary, we obtain~that
\Bsist
  && 2\eta \intQt |D\vvk|^2
  + \nu \intQt |\vvk|^2
  \non
  \\
  && = 2\eta \intQt D\vvk:\nabla\uuk
  + \nu \intQt \vvk \cdot \uuk
  \non
  \\
  && \quad {}
  + \intQt \muk \nabla\phik \cdot \vvk
  + \intQt \Ns(\phik,\sigk) \nabla\sigk \cdot \vvk
  \non
  \\
  && \quad {}
  - \intQt \muk \nabla\phik \cdot \uuk
  - \intQt \Ns(\phik,\sigk) \nabla\sigk \cdot \uuk
  + \intQ \heps \nabla\phik \cdot (\vvk-\uuk) \,.
  \qquad
  \label{testprimak}
\Esist
At the same time, we test \eqref{secondak} and \eqref{terzak} 
by $\muk$ and $\dt\phik$, respectively, and integrate with respect to time.
We obtain~that
\Bsist
  && \frac\eps 2 \iO |\muk(t)|^2
  + \intQt |\nabla\muk|^2 
  \non
  \\
  && = - \intQt \dt\phik \, \muk
  + \intQt \Sph(\phik,\sigk) \muk
  - \intQt \nabla\phik \cdot \vvk \, \muk
  \an{-} \intQt \an{\phi_k}\geps\muk\,,
  \label{testsecondak}
\Esist
as well as
\begin{align}
  & \eps \intQt |\dt\phik|^2
  + \frac 12 \iO |\nabla\phik(t)|^2 
  + \iO \Feps(\phik(t)) 
  \pier{{}+ \intQt \Nph(\phik,\sigk) \dt\phik{}}  
  \non
  \\
  & = \frac 12 \iO |\nabla\phik(0)|^2  
  + \iO \Feps(\phik(0)) 
  + \intQt \muk \dt\phik \an{{}+ \intQt \heps \dt\phik} \,.
  \label{testterzak}
\end{align}
By recalling that $\Ns(\phik,\sigk)=\sigk+\chi(1-\phik)$ by \eqref{idnut}, 
testing \eqref{quartak} by $\Ns(\phik,\sigk)$
and integrating with respect to time, we \an{also} obtain~that
\begin{align}
  & \intQt \dt\sigk \, \Ns(\phik,\sigk)
  + \intQt |\nabla\Ns(\phik,\sigk)|^2
  + \kappa \intSt |\sigk|^2
  \non
  \\
  & = \intQt \Ss (\phik,\sigk) \Ns(\phik,\sigk)
  - \intQt \nabla\sigk \cdot \vvk \, \Ns(\phik,\sigk)
  - \intQt \sigk \geps \Ns(\phik,\sigk)
  \non
  \\
  & \quad {}
  + \kappa \intSt \sig_\Sigma \bigl( \sigk + \chi(1-\phik) \bigr)
  - \kappa \chi \intSt \sigk (1-\phik)\an{.}
  \label{testquartak}
\end{align}
At this point, we add \accorpa{testprimak}{testquartak} to each other
and notice that some cancellations occur.
Moreover, we \pier{combine the term \last{involving $N_\phi(\phik,\sig_k)$ from} \eqref{testterzak}  
with the first one of~\eqref{testquartak}.} 
By recalling \eqref{idnut}, we have~that
\Beq
  \intQt \Nph(\phik,\sigk) \dt\phik
  + \intQt \dt\sigk \, \Ns(\phik,\sigk)
  = \iO N(\phik(t),\sigk(t))
  - \iO N(\phik(0),\sigk(0)).
  \non
\Eeq
Next, we observe that the Korn inequality \eqref{korn} implies~that
\Beq
  2\eta \intQt |D\vvk|^2
  + \nu \intQt |\vvk|^2
  \geq \alpha \iot \normaV{\vvk(s)}^2 \, ds\,,
  \quad \hbox{where} \quad
  \alpha := \min\{2\eta,\nu\} / C_K \,.
  \non
\Eeq
Finally, as for the term involving~$\Feps$, we apply \eqref{coerceps} with $M=1$,
and obtain~that
\Beq
  \iO \Feps(\phik(t))
  \geq \iO |\phik(t)|^2
  - c \,,
  \non
\Eeq
provided that $\eps$ is small enough (as in the lemma)\an{: from} now on, it is understood that $\eps$ satisfies this smallness condition.
Regarding the \rhs,
we start by treating the nontrivial terms coming from the identity~\eqref{testprimak}. 
The symbol $\delta$ denotes a positive parameter whose value is chosen later on.
By recalling \eqref{disuguk}, we have~that
\Bsist
  && 2\eta \intQt D\vvk:\nabla\uuk
  + \nu \intQt \vvk \cdot \uuk  
  \non
  \\
  && \leq \delta \iot \normaV{\vvk(s)}^2 \, ds
  + c_\delta \iot \normaV{\uuk(s)}^2 \, ds
  \leq \delta \iot \normaV{\vvk(s)}^2 \, ds
  + c_\delta \,.
  \non
\Esist
Next, \pier{with the help of} the \Holder, Sobolev, Poincar\'e, and Young inequalities,
\pier{we obtain~that}
\Bsist
  && - \intQt \muk \nabla\phik \cdot \uuk
  \leq \iot \norma{\muk(s)}_4 \, \an{\norma{\nabla\phik(s)}} \, \norma{\uuk(s)}_4 \, ds
  \non
  \\
  && \leq \delta \intQt |\nabla\muk|^2
  + c_\delta \iot \an{\norma{\nabla\phik(s)}^2 \, \norma{\uuk(s)}_4^2} \, ds
  \leq \delta \intQt |\nabla\muk|^2
  + c_\delta \intQt |\nabla\phik|^2 \,.
  \non
\Esist
For the next term, we also apply the compactness inequality \eqref{compact} to $\Ns(\phik,\sigk)$ and have~that
\Bsist
  && - \intQt \Ns(\phik,\sigk) \nabla\sigk \cdot \uuk
  = - \intQt \Ns(\phik,\sigk) \bigl( \nabla\Ns(\phik,\sigk) + \chi \nabla\phik \bigr) \cdot \uuk
  \non
  \\
  && \leq \iot \norma{\Ns(\phik(s),\sigk(s))}_4 \,
    \bigl( \an{\norma{\nabla\Ns(\phik(s),\sigk(s))}} + \chi \an{\norma{\nabla\phik(s)}}
    \bigr) \norma{\uuk(s)}_4 \, ds
  \non
  \\
  && \leq \delta \intQt \bigl( |\nabla\Ns(\phik,\sigk)|^2 + |\nabla\phik|^2 \bigr)
  + c_\delta \intQt |\Ns(\phik,\sigk)|^2 \,.
  \non
\Esist
Similarly, using \eqref{hpheps} and the \pier{Sobolev inequality \eqref{sobolev} for 
$\heps$ and $\vvk$}, we have that
\Beq
  \an{\intQt} \heps \nabla\phik \cdot (\vvk-\uuk) 
  \leq \delta \pier{{}\iot \norma{\vvk(s)}_{\VV}^2 \, ds{}}
  + c_\delta \iot \normaV{\heps(s)}^2 \, \norma{\nabla\phik(s)}^2 \, ds
  + c_\delta \,.
  \non
\Eeq
We notice that the $L^1$ norm of the function $s\mapsto\normaV{\heps(s)}^2$ 
is bounded by a constant independent of $\eps$ by~\eqref{hpheps}.
Now, among the volume integrals that come from \accorpa{testsecondak}{testquartak} and should be estimated,
just the last one \juerg{on the \rhs\ of}  \eqref{testterzak} needs some treatment.
Indeed, all the others can easily be dealt with by \pier{virtue of} the Young inequality, 
possibly combined with other estimates, like \last{\eqref{hpSbis} or \eqref{poincare}}, without any difficulty.
We have~that
\Bsist
  && \intQt \heps \dt\phik
  = - \intQt \dt\heps \, \phik
  + \iO \heps(t) \phik(t)
  - \iO \heps(0) \phik(0)
  \non
  \\
  && \leq \intQt ( |\dt\heps|^2 + |\phik|^2 )
  + \delta \, \iO |\phik(t)|^2 + c_\delta \iO |\heps(t)|^2
  + \iO |\an{\phi_k}(0)|^2 + \iO |\heps(0)|^2
  \non
  \\
  && \leq \intQt |\phik|^2
  + \delta \, \iO |\phik(t)|^2
  + \iO |\an{\phi_k}(0)|^2
  + c_\delta\,,
  \non
\Esist
where we owe  to \eqref{hpheps} for the last inequality.
Now, we \an{move} to the surface integrals.
We have~that
\Bsist
  && \kappa \intSt \sig_\Sigma \bigl( \sigk + \chi(1-\phik) \bigr)
  - \kappa \chi \intSt \sigk (1-\phik)
  \leq \frac \kappa 2 \intSt |\sigk|^2
  + c \intSt |\phik|^2 
  + c
  \non
  \\
  && \leq \pier{\frac \kappa 2} \intSt |\sigk|^2
  + \pier{c} \iot \normaV{\phik(s)}^2 \, ds
  + \pier{c}
  \non
  \\
  && \leq \pier{\frac \kappa 2} \intSt |\sigk|^2
  + \pier{c} \intQt |\phik|^2 
  + \pier{c} \intQt |\nabla\phik|^2
  + \pier{c} \,.
  \non
\Esist
Finally, all of the terms involving the initial values
can easily be estimated by accounting for~\eqref{projk}
and recalling \eqref{disugBetaeps} and \eqref{defN} 
to treat the convex part $\Betaeps$ of $\Feps$ and the term involving~$N$.
At this point, by collecting everything, choosing $\delta$ small enough, and applying the Gronwall lemma,
we obtain that
\Bsist
  && \norma\vvk_{L^2(0,T_k;\VV)}
  + \norma{\nabla\muk}_{L^2(0,T_k;\an{\HH})}
  + \norma\phik_{L^\infty(0,T_k;V)}
  + \norma{\Feps(\phik)}_{L^\infty(0,T_k;\Luno)}
  \non
  \\[1mm]
  && \quad {}
  + \norma{N(\phik,\sigk)}_{L^\infty(0,T_k;H)}
  + \norma{\Ns(\phik,\sigk)}_{L^2(0,T_k;\pier{V})}
  \non
  \\[1mm]
  && \quad {}
  + \eps^{1/2} \, \norma\muk_{L^\infty(0,T_k;H)}
  + \eps^{1/2} \, \norma{\dt\phik}_{L^2(0,T_k;H)}
  \leq c \,.
  \non
\Esist
In particular, this proves \eqref{fromprimastima}, so that $T_k=T$.
Then, by recalling \Nutrient\ and the Poincar\'e inequality once more and rearranging,
we conclude~that
\Bsist
  && \norma\vvk_{\L2\VV}
  + \norma\muk_{\L2\Vz}
  + \norma\phik_{\L\infty V}
  + \norma\sigk_{\L\infty H\cap\L2V}
  \non
  \\
  && \quad {}
  + \eps^{1/2} \, \norma\muk_{\L\infty H}
  + \eps^{1/2} \, \norma{\dt\phik}_{\L2H}
  \leq c \,.
  \qquad
  \label{primastima}
\Esist

\step
Second a priori estimate

We now aim at recovering \an{an} estimate for the pressure~$\pk$. 
\an{Thus, w}e construct $\qqk\in\L2\VV$ such that
\Bsist
  && \div\qqk(t) = \pk(t)
  \quad \hbox{in $\Omega$}
  \aand
  \qqk(t)\suG = \frac 1{|\Gamma|} \bigl( \textstyle\iO \pk(t) \bigr) \nn\,,
  \non
  \\
  && \normaV{\qq(t)} \leq C \, \norma{\pk(t)}\,,
  \non
\Esist
\aat\ and some constant~$C\an{{}>0}$.
To this end, it suffices to apply Lemma~\ref{Lemmau} with an obvious choice of $f$ and~$\aa$.
Then, we test \eqref{primak}, written at the time $t$, by~$\qqk(t)$.
However, we avoid writing the time $t$ for brevity.
\pier{Recalling \eqref{def:T} and \eqref{eqdivk}, 
we} obtain~that
\Bsist
  &&\iO \bigl(
    2\eta D\vvk:\nabla\qqk
    + \lambda \geps \, \div\qqk
    \pier{{}-{}} \pk \div\qqk
    \pier{{} + \nu \hskip1pt\vvk \cdot \qqk} 
  \bigr)
  \non
  \\
  &&= \iO \bigl(
    \an{(\muk + \heps) \nabla\phik}
    + \Ns(\phik,\sigk) \nabla\sigk
  \bigr) \cdot \qqk\,,
  \non
\Esist
and the definition of $\qqk$\juerg{, as well as} the \Holder, Sobolev, and Young inequalities, yield
\begin{align*}
   \norma\pk^2 
  &\leq c \bigl( \normaV\vvk + \norma\geps \bigr) \normaV\qqk
  \non
  \\
  & \quad {}
  + (\norma\muk_4 + \norma\heps_4) \, \norma{\nabla\phik} \, \norma\qqk_4
  + \norma{\Ns(\phik,\sigk)}_3 \, \norma{\nabla\sigk} \, \norma\qqk_6
  \non
  \\
  &\leq c \bigl( \normaV\vvk + \an{\norma\geps} \bigr) \norma\pk
  \non
  \\
  & \quad {}
    + c (\norma\muk_4 + \norma\heps_4) \, \normaV\phik \, \norma\pk
   \an{{} + c \, \norma{\Ns(\phik,\sigk)}_3 \, \normaV\sigk \, \norma\pk}
  \non
  \\
  & \leq \frac 12 \, \norma\pk^2
  + c \bigl( \normaV\vvk^2 + \an{\norma\geps}^2 \bigr)
  \non
  \\
  & \quad {}  
   + c (\norma\muk^2_{\pier{\Vz}}  + \norma\heps^2_{\pier{V}} ) \normaV\phik^2
    + c \, \norma{\Ns(\phik,\sigk)}_3^2 \, \normaV\sigk^2 \,.
  \non
\end{align*}
Now, we rearrange, take the power of exponent $2/3$, and integrate over $(0,T)$. 
By accounting for \eqref{primastima} and the Young inequality once more, we deduce~that
\begin{align}
  & \ioT \an{\norma\pk^{4/3}} \, dt
  \non
  \\
  & \leq c \ioT \bigl(
    \normaV\vvk^{4/3}
    + \pier{\norma\geps^{4/3}}
    + (\pier{\norma\muk_{\Vz}^{4/3} + \norma\heps_V^{4/3}}) \normaV\phik^{4/3}
    + \norma{\Ns(\phik,\sigk)}_3^{4/3} \, \normaV\sigk^{4/3}
  \bigr) \, dt
  \non
  \\
  & \leq c \ioT \norma{\Ns(\phik,\sigk)}_3^{4/3} \, \normaV\sigk^{4/3} \, dt
  + c
  \leq c \ioT \bigl( \norma{\Ns(\phik,\sigk)}_3^4 + \normaV\sigk^2 \bigr) \, dt
  + c
  \non
  \\
  & \leq c \ioT \norma{\Ns(\phik,\sigk)}_3^4 \, dt
  + c \,,
  \non
\end{align}
and it remains to estimate the last integral.
By interpolation, we have the continuous embedding
\Beq
  \L\infty H \cap \L2{\Lx6} \emb \L4{\Lx3} \,.
  \non
\Eeq
Since \eqref{primastima} implies that $\Ns(\phik,\sigk)$ is bounded in $\an{\L\infty H}\cap\L2V$
and the continuous embedding $V\emb\Lx6$ holds, the integral at hand is uniformly bounded.
Therefore, we have proved~that
\Beq
  \norma\pk_{\L{4/3}H} \leq c \,.
  \label{secondastima}
\Eeq
By the way, the argument used for $\Ns(\phik,\sigk)$ also applies to~$\sigk$, so~that
\Beq
  \norma\sigk_{\L4{\Lx3}} \leq c \,.
  \label{stimasigk}
\Eeq

\step
Third a priori estimate

We test \eqref{terzak} by the admissible function $-\Delta\phik$
and integrate in time.
We obtain~that
\Bsist
  && \frac \eps 2 \iO |\nabla\phik(t)|^2
  + \intQt |\Delta\phik|^2
  + \intQt \betaeps'(\phik) |\nabla\phik|^2
  \non
  \\
  && = \frac \eps 2 \iO |\nabla\phik\pier{(0)}|^2
  + \intQt f (-\Delta\phik)
  \quad \hbox{where} \quad
  f := \muk \an{{}+ \heps{}} - \Nph(\phik,\sigk)  - \pi(\phik) \,.
  \non
\Esist
As for the first term on the \rhs, we recall \eqref{hpzero} and Remark~\ref{Projk}, 
to see that it is bounded.
Since $f$ is bounded in $\L2H$ by \eqref{primastima},
we deduce that the same holds for $\Delta\phik$.
Then, elliptic regularity yields that
\Beq
  \norma\phik_{\L2W} \leq c \,.
  \label{terzastima}
\Eeq

\step
Fourth a priori estimate

We take any $\zeta\in\L4V$ and define $\zetak$
\juerg{as follows}: $\an{\zetak}(t)$ is the $H$-projection of $\zeta(t)$ onto~$\Vk$ \aat.
Then, \aat, we test \eqref{quartak}, written at the time $t$, by $\zetak(t)$. 
However, we do not write the time~$t$ for simplicity.
By also accounting for \eqref{eqdivk}, we have~that
\Bsist
  && \iO \dt\sigk \, \zetak
  = - \iO \nabla\Ns(\phik,\sigk) \cdot \nabla\zetak
  + \iO \Ss(\phik,\sigk) \zetak
  \non
  \\
  && \quad {}
  - \iO\bigl( \nabla\sigk \cdot \vvk + \sigk \div\vvk \bigl) \zetak
  \an{{}+{}} \kappa \iG (\an{\sig_\Sigma - \sigk }) \zetak \,.
  \non
\Esist
Since $\dt\sigk(t)\in\Vk$ and $\zetak(t)$ coincides with the $V$-projection of~$\zeta(t)$ 
as explained in Remark~\ref{Projk}, 
we can replace $\zetak$ by~$\zeta$ (still omitting the time) on the \lhs\ 
and obtain 
\Beq
  \iO \dt\sigk \, \zetak
  = \iO \dt\sigk \, \zeta \,.
  \non
\Eeq
By the same remark, we have that $\normaV\zetak\leq\normaV\zeta$ (see \eqref{genprojk}).
Hence, we infer that
\Beq
  - \iO \nabla\Ns(\phik,\sigk) \cdot \nabla\zetak
  + \iO \Ss(\phik,\sigk) \zetak
  \leq c \bigl( \normaV\phik + \normaV\sigk + 1 \bigr) \normaV\zeta \,.
  \non
\Eeq
To deal with the next term, we integrate by parts and obtain
\Bsist
  && - \iO \bigl( \nabla\sigk \cdot \vvk + \sigk \div\vvk \bigl) \zetak
  = - \iO (\div(\sigk\vvk)) \, \zetak
  \non
  \\
  && = \iO \sigk \, \vvk \cdot \nabla\zetak
  - \iG \sigk \, \zetak \, \vvk\cdot\nn
  \non
  \\[2mm]
  && \leq \norma\sigk_3 \, \norma\vvk_6 \, \norma{\nabla\zetak}
  + \norma{\sigk\suG} \, \norma{\vvk\cdot\nn}_4 \, \norma{\zetak\suG}_4 \,.
  \label{divergence}
\Esist
We can replace the $L^6$ norm by the $V$ norm since $\an{ V \emb \Lx6}$,
and the $L^2$ norm of $\nabla\zetak$ by~$\normaV\zeta$.
Moreover, since the two-dimensional embedding $\HxG{1/2}\emb\LxG4$ holds true
and the trace operator is continuous from $V$ to $\HxG{1/2}$,
we can estimate the last product as follows:
\Beq
  \norma{\sigk\suG} \, \norma{\vvk\cdot\nn}_4 \, \norma{\zetak\suG}_4
  \leq c \, \norma{\sigk\suG} \, \normaV\vvk \, \normaV\zeta \,.
  \label{divergencebis}
\Eeq
Similarly, we have that
\Beq
  \kappa \iG (\an{\sig_\Sigma - \sigk }) \zetak 
  \leq c \bigl(\an{\norma{\sig_\Sigma}  + \norma{\sigk\suG} } \bigr) \, \normaV\zetak 
  \leq c \bigl( \an{\norma{\sig_\Sigma} + \normaV\sigk} \bigr) \, \normaV\zeta \,,
  \non
\Eeq
where, for clarity, we point out that $\norma{\sig_\Sigma}$ here means the norm of $\sig_\Sigma(t)$ in $\LdueG$.
At this point, we collect all these \last{in}equalities and estimates and integrate over~$(0,T)$.
Omitting the integration variable $t$ for brevity, we have that
\Bsist
  && \intQ \dt\sigk \, \zeta
  \leq c \ioT \bigl( \normaV\phik + \normaV\sigk + 1 \bigr) \normaV\zeta \, dt
  + c \ioT \norma\sigk_3 \, \normaV\vvk \, \normaV\zeta \, dt
  \non
  \\
  && \quad {}
  + c \ioT \norma{\sigk\suG} \, \normaV\vvk \, \normaV\zeta \, dt
  + c \ioT \bigl( \an{\norma{\sig_\Sigma} + \normaV\sigk } \bigr) \, \normaV\zeta \, dt \,.
  \non
\Esist
Therefore, by using the \Holder\ inequality, we have~that
\Bsist
  && \intQ \dt\sigk \, \zeta
  \leq c \bigl( \norma\phik_{\L2V} + \norma\sigk_{\L2V} \an{+1}\bigr) \norma\zeta_{\L2V}
  \non
  \\
  && \quad {}
  + c \, \norma\sigk_{\L4{\Lx3}} \, \norma\vvk_{\L2\VV} \, \norma\zeta_{\L4V}
  \non
  \\
  && \quad {}
  + c \, \norma{\sigk\suG}_{\L4\LdueG} \, \norma{\an{\vvk}}_{\L2\VV} \, \norma\zeta_{\L4V}
  \non
  \\
  && \quad {}
  + c \, \bigl( \an{\norma{\sig_\Sigma}_{\LS2}  + \norma\sigk_{\L2V} }\bigr) \norma\zeta_{\L2V} \,.
  \non
\Esist
Finally, we account for \eqref{primastima}, \eqref{stimasigk}, and Lemma~\ref{Traccia},
\an{to} conclude that
\Beq
  \intQ \dt\sigk \, \zeta
  \leq c \, \norma\zeta_{\L4V} 
  \quad \hbox{for every $\zeta\in\L4V$}.
  \non
\Eeq
This means that $\dt\sigk$ is bounded in the dual space of~$\L4V$, i.e., that
\Beq
  \norma{\dt\sigk}_{\L{4/3}\Vp} \leq c \,.
  \label{quartastima}
\Eeq

\step
Fifth a priori estimate

We aim at estimating \an{the time derivative} $\dt(\eps\muk+\phik)$.
To this end, we take any $\badphi\in\L2\Vz$
and define $\badphik:Q\to\erre$ by setting \aat:
$\badphik(t)$ is the $H$-projection of $\badphi(t)$ onto~$\Vzk$.
Then, we rearrange \eqref{secondak}, written at the time~$t$, 
test it by~$\badphik(t)$, and integrate over~$(0,T)$.
We obtain~that
\Bsist
  && \intQ \dt(\eps\muk+\phik) \, \badphik
  = - \intQ \nabla\muk \cdot \nabla\badphik
  + \intQ \Sph(\phik,\sigk) \badphik
  - \intQ (\nabla\phik \cdot \vvk + \an{\phik}\geps) \badphik \,.
  \non
\Esist
\juerg{For the first two terms on the \rhs, we invoke the \Lip\ continuity of $\Sph$ and \eqref{primastima}
to immediately obtain} that
\Bsist
  && - \intQ \nabla\muk \cdot \nabla\badphik
  + \intQ \Sph(\phik,\sigk) \badphik
  \non
  \\[2mm]
  && \leq c \, \norma\muk_{\L2\Vz} \, \norma\badphik_{\L2\Vz}
  + c \bigl( \norma\phik_{\L2H} + \norma\sigk_{\L2H} + 1 \bigr) \norma\badphik_{\L2H}
  \non
  \\[2mm]
  && \leq c \, \norma\badphik_{\L2\Vz},
  \non
\Esist
while the last one needs some care.
By the \Holder\ and Sobolev inequalities, \eqref{hpgeps} and \eqref{primastima} once more, we obtain~that
\Bsist
  && - \intQ (\nabla\phik \cdot \vvk + \phik\geps) \badphik
  \non
  \\
  && \leq \ioT \norma{\nabla\phik(t)} \, \norma{\vvk(t)}_4 \, \norma{\badphik(t)}_4 \, dt
  + \norma\phik_{\L2 {\last{H}}} \, \norma\geps_\infty \, \norma\badphik_{\L2 {\last{H}}}
  \non
  \\
  && \leq \norma\phik_{\L\infty V} \, \norma \vvk_{\L2\VV} \, \norma\badphik_{\L2\Vz} 
  + c \, \norma\phik_{\L2{\last{H}}} \, \norma\badphik_{\L2{\last{H}}}
  \non
  \\
  && \leq c \, \norma\badphik_{\L2\Vz} \,.
  \non
\Esist
Since $\norma\badphik_{\L2\Vz}\leq\norma\badphi_{\L2\Vz}$ by Remark~\ref{Projk}
and we can replace $\badphi$ by~$-\badphi$, we infer~that
\Beq
  \Bigl| \intQ \dt(\eps\muk+\phik) \, \badphik \Bigr|
  \leq c \, \norma\badphi_{\L2\Vz} \,.
  \label{perquintastima}
\Eeq
Unfortunately, just $\dt\muk$ is $\Vzk$-valued while $\dt\phik$ is not, 
so that we only have~that
\Beq
  \intQ \dt(\eps\muk+\phik) \, \badphi
  = \intQ \dt(\eps\muk+\phik) \, \badphik
  + \intQ \dt\phik (\badphi-\badphik).
  \non
\Eeq
However, on account of \eqref{primastima} and Remark~\ref{Projk} once more, we can write
\Beq
  \Bigl| \intQ \dt\phik (\badphi-\badphik) \Bigr|
  \leq \norma{\dt\phik}_{\L2H} \bigl( \norma\badphi_{\L2H} + \norma\badphik_{\L2H} \bigr)
  \leq c_\eps \, \norma\badphi_{\L2\Vz} \,.
  \non
\Eeq
Combining this with \eqref{perquintastima} yields
\Beq
  \norma{\dt(\eps\muk+\phik)}_{\L2\Vzp} \leq c_\eps \,.
  \label{quintastima}
\Eeq

\step
\juerg{Passage to the limit as $k \to \infty$}

We collect the basic estimates we have proved, namely,
\eqref{primastima}, \eqref{secondastima}, \eqref{terzastima}, \eqref{quartastima}, and \eqref{quintastima},
and apply well-known weak and weak star compactness results.
Since $\eps$ is fixed, for some (not relabeled) subsequence and suitable limit functions, we have\an{, as $k\to\infty$,}
\begin{align}
  & \vvk \to \vveps
  \quad \hbox{weakly in $\L2\VV\emb\L2{\LLx4}$}\,,
  \label{convvvk}
  \\
  & \pk \to \peps
  \quad \hbox{weakly in $\L{4/3}H$}\,,
  \label{convpk}
  \\
  &\muk \to \mueps
  \quad \hbox{weakly star in $\L\infty H\cap\L2\Vz$}\,,
  \label{convmuk}
  \\
  &\phik \to \phieps
  \quad \hbox{weakly star in $\H1H\cap\L\infty V\cap\L2W$}\,,
  \label{convphik}
  \\
  &\sigk \to \sigeps
  \quad \hbox{weakly star in $\W{1,4/3}\Vp\cap\L\infty H\cap\L2V$}\,,
  \qquad
  \label{convsigk}
  \\
  &\eps\muk+\phik \to \eps\mueps+\phieps
  \quad \hbox{weakly in $\H1\Vzp$} \,.
  \label{convcombk}
\end{align}
\pcol{In view of \eqref{convphik}, note that \eqref{convcombk}
is equivalent to $\muk \to \mueps$ weakly in $\H1\Vzp$.}
Then, the quintuple $\soluzeps$ satisfies \Regsoluzeps, 
as well as the initial conditions~\pier{\eqref{cauchyeps}}
(recall the Remarks~\ref{Regmueps} and~\ref{Projk}),
and we aim at proving that it solves the whole approximating problem.
By linearity (cf.~\eqref{idnut}), we deduce that\an{, as $k \to \infty$,}
\Bsist
  && \Nph(\phik,\sigk) \to \Nph(\phieps,\sigeps)
  \aand
  \Ns(\phik,\sigk) \to \Ns(\phieps,\sigeps)
  \non
  \\
  && \quad \hbox{weakly star in $\W{1,4/3}\Vp\cap\L\infty H\cap\L2V$} \,.
\non  
\Esist
Moreover, by strong compactness
(see, e.g., \cite[Thm.~5.1, p.~58]{Lions} and \cite[Sect.~8, Cor.~4]{Simon})
we \juerg{may without loss of generality assume}~that
\pier{%
\begin{align}
& \phik \to \phieps
  \quad \hbox{strongly in $\L2{\Wx{1,4}}$ and \aeQ}\,,
  \label{strongphik}
  \\
  & \sigk \to \sigeps
  \quad\hbox{strongly in $\L2{\Lx4}$ and \aeQ}\,,
  \label{strongsigk}
  \\
  & \eps\muk+\phik \to \eps\mueps+\phieps
  \quad \hbox{strongly in $\L2{\Lx4}$} \,.
   \label{strongcombk}
\end{align}
}%
Notice that combining \eqref{convvvk} and \eqref{convmuk} with \eqref{strongphik} yields that
\pier{%
\begin{align}
  & \nabla\phik \cdot \vvk  \to \nabla\phieps \cdot \vveps
  \quad \hbox{weakly in $\L1{\an{H}}$}\,, 
    \label{weakprod:1}
   \\
	& \muk \nabla\phik  \to \mueps \nabla\phieps
  \quad \hbox{weakly in $\L1{\an{\HH}}$} \,.
  \label{weakprod:2}
\end{align}
}%
On the other hand, the \pier{above convergence properties and the \Lip\ continuity imply that}\an{, as $k\to\infty$,}
\Beq
  \Ss(\phik,\sigk) \to \Ss(\phieps,\sigeps)
  \aand
  \Feps'(\phik) \to \Feps'(\phieps)\,,
  \quad \hbox{strongly in $\L2{\Lx4}$}.
  \label{strongnonlin}
\Eeq
Now, we \an{claim}~that
\Beq
  \zetak \to \zeta
  \enskip \hbox{strongly in $\L\infty V$}
  \quad \hbox{implies that} \quad
  \intQ (\nabla\sigk \cdot \vvk) \zetak
  \to \intQ (\nabla\sigeps \cdot \vveps) \zeta \,.
  \label{convprod}
\Eeq
Assume $\zetak$ and $\zeta$ as said.
From \eqref{convsigk} and the strong compactness results already quoted,
we deduce that\an{, as $k\to\infty$,}
\Beq
  \sigk \to \sigeps
  \quad \hbox{strongly in $\L2{\Hx{3/4}}$} .
  \non
\Eeq
From this, we deduce that\an{, as $k\to\infty$,}
\Beq
  \nabla\sigk \to \nabla\sigeps
  \quad \hbox{strongly in $\L2{\HHx{-1/2}}$}\,, 
  \label{strongnablasigma}
\Eeq
\pier{owing to the} interpolation theory in Hilbert spaces (see, e.g., \cite[Ch.~I]{LioMag}).
We know that $\nabla$~is continuous from $\an{V}$ into $\an{\HH}$ 
and from $\an{H}$ into $\HHx{-1}=(\HH^1_0(\Omega))^*$.
Then, it is continuous 
from $\Hx{3/4}=(\an{V},\an{H})_{1/4}$ into $(\an{\HH},(\HH^1_0(\Omega))^*)_{1/4}$.
Since $1/4<1/2$, the last space is $(\HH^{1/4}_0(\Omega))^*=(\HHx{1/4})^*$,
\juerg{which} is embedded in $\HHx{-1/2}$ (see~\accorpa{embeddings}{dualeHs}).
Then, \eqref{strongnablasigma} follows.
On the other hand, we have that 
\Bsist
  && \norma{\nabla(\zetak\vvk)}_{\L2{\LLx{3/2}}}
  \non
  \\
  && \leq \norma{\nabla\zetak}_{\L\infty{\an{\HH}}} \, \norma\vvk_{\L2{\LLx6}}
  + \norma\zetak_{\L\infty{\Lx6}} \, \norma{\nabla\vvk}_{\L2{\an{\HH}}}
  \non
  \\
  && \leq c \, \norma\zetak_{\L\infty{\an{V}}} \, \norma\vvk_{\L2\VV}
  \leq c \,,
  \non
\Esist
so that $\{\zetak\vvk\}$ has a weak limit in $\L2{\WWx{1,3/2}}$.
Since $\zetak\vvk\to\zeta\vveps$ weakly in $\an{\L2 {\HH}}$, we deduce that
\Beq
  \zetak\vvk \to \zeta\vveps
  \quad \hbox{weakly in $\L2{\WWx{1,3/2}}\emb\L2{\HHx{1/2}}$}\,,
  \non
\Eeq
where the last (three-dimensional) embedding follows from a particular case of the embedding
\Beq
  \Wx{1,q} \emb \Hx s = \Wx{s,2},
  \quad \hbox{where} \quad
  s \in (0,1) , \quad
  q \in [\an{1},+\infty),
  \aand
  s - \frac 32 \leq 1 - \frac 3q \,.
  \non
\Eeq
Therefore, \eqref{convprod} holds true since $\L2{\HHx{-1/2}}=(\L2{\HHx{1/2}})^*$.
By the same argument (by~replacing the velocities by the $\Ns$ terms, essentially),
one obtains~that
\Beq
  \intQ  \an{\Ns(\phik,\sigk)\nabla\sigk\cdot \bzeta
  \to \intQ \Ns(\phieps,\sigeps)\nabla\sigeps \cdot \bzeta}
  \quad \hbox{for every $\bzeta\in\L\infty\VV$} \,.
  \label{convprodbis}
\Eeq

\step
Conclusion of the proof of Theorem~\ref{Existenceeps}

At this point, we can verify that the \an{quintuple} just found solves the equations of problem \Pbleps.
Clearly, \eqref{eqdiveps} follows from \eqref{eqdivk} \pier{and \eqref{convvvk}.}
As for \eqref{primaeps}, it suffices to check that some equivalent formulation is satisfied.
This is the case if we take
\Bsist
  && \intQ \T(\vveps,\peps) : \nabla\bzeta 
  + \nu \intQ \vveps \cdot \bzeta 
  = \intQ (\mueps+\heps) \nabla\phieps \cdot \bzeta
  + \intQ \Ns(\phieps,\sigeps) \nabla\sigeps \cdot \bzeta
  \non
  \\
  && \quad \hbox{for every $\bzeta\in\L\infty\VV$} .
  \label{intprimaeps}
\Esist
Indeed, this equation actually is satisfied due to the convergence properties just mentioned.
Similarly, instead of considering \accorpa{secondaeps}{quartaeps},
we deal with some time integrated version of theirs (analogous to \eqref{intprimaeps}),
but with step \pier{functions} as test functions 
(this is convenient at~least for some of the equations, and for simplicity we choose step functions for all of them).
So,~we assume that $\badphi$ is a $\Vz$-valued step function and that
$z$ and $\zeta$ are $V$-valued step functions.
However, we have to be careful, since we \an{are starting} from the discrete setting.
So, given~$\badphi$, $z$ and $\zeta$ as said, we introduce $\badphik$ 
as we did in proving our fifth estimate, i.e.,
by \juerg{defining
$\badphik(t)$ as} the $H$-projection of $\badphi(t)$ onto~$\Vzk$, \aat.
Similarly, we define $\zk$ and $\zetak$ starting from $z$ and~$\zeta$,
now with $\Vk$ instead of~$\Vzk$.
By accounting for Remark~\ref{Projk},
we point out at once~that\an{, as $k\to\infty$,}
\Beq
  \badphik \to \badphi , \quad
  \zk \to z,
  \aand
  \zetak \to \zeta,  \quad \hbox{strongly in $\L\infty V$},
  \label{convtest}
\Eeq
since $\badphi$, $z$ and $\zeta$ have a finite number of values.
\an{Next}, we test \eqref{secondak}, written at the time $t$, by $\badphik(t)$ and integrate over~$(0,T)$.
\an{Upon} rearranging, we obtain~that
\Beq
  \intQ \dt(\eps\muk+\phik) \, \badphik
  + \intQ \nabla\muk \cdot \nabla\badphik
  = \intQ \Sph(\phik,\sigk) \badphik
  - \intQ \bigl( \nabla\phik \cdot \vvk + \phik \geps \bigr) \badphik \,.
  \non
\Eeq
On account of \an{the strong convergence in} \eqref{convtest}, \pier{recalling \eqref{convcombk} and \eqref{weakprod:1}, and}
letting $k$ tend to infinity, we conclude~that
\Bsist
  && \ioT \< \dt(\eps\mueps+\phieps)(t) , \badphi(t) >_{\Vz}
  + \intQ \nabla\mueps \cdot \nabla\badphi
  \non
  \\
  && = \intQ \Sph(\phieps,\sigeps) \badphi
  - \intQ \bigl( \nabla\phieps \cdot \vveps + \phieps \geps \bigr) \badphi \,.
  \qquad
  \label{intsecondaeps}
\Esist
What we have obtained is an integrated version of \eqref{secondaeps},
which is equivalent to \eqref{secondaeps} itself
since it holds for every $\Vz$-valued step function~$\badphi$.
Similarly, we test \eqref{terzak} and \eqref{quartak}
written at the time $t$ by $\zk(t)$ and~$\zetak(t)$, respectively, 
and integrate over~$(0,T)$.
We obtain~that
\Bsist
	&& \eps \intQ \dt\phik \, \zk 
    + \intQ \nabla\phik \cdot \nabla \zk 
	+ \intQ \Feps'(\phik) \, \zk 
	= \intQ \bigl( \muk \an{{}+ \heps{} } - \Nph(\phik,\sigk) \bigr) \zk 
	\non
	\\
	&& \intQ \dt \sigk \, \zetak
	+ \intQ \nabla\Ns(\phik,\sigk) \cdot \nabla\zetak
	\non
	\\
	&& = \intQ \Ss (\phik,\sigk) \zetak
	- \iO (\nabla\sigk \cdot \vvk + \sigk \geps) \pier{\zetak}
	+ \kappa \intS (\sig_\Sigma - \sigk) \zetak \,. 
	\non
\Esist
By accounting for \last{\eqref{convprod} and \eqref{convtest}}, we conclude~that\an{, as $k\to\infty$,}
\Bsist
	&& \eps \intQ \dt\phieps \, z
    + \intQ \nabla\phieps \cdot \nabla z
	+ \intQ \Feps'(\phieps) \, z
	= \intQ \bigl( \mueps \an{{}+ \heps{} } - \Nph(\phieps,\sigeps) \bigr) z\,,
	\qquad
	\label{intterzaeps}
	\\
	&& \intQ \dt \sigeps \, \zeta
	+ \intQ \nabla\Ns(\phieps,\sigeps) \cdot \nabla\zeta
	\non
	\\
	&& = \intQ \Ss (\phieps,\sigeps) \zeta
	- \iO (\nabla\sigeps \cdot \vveps + \sigeps \geps) \zeta
	+ \kappa \intS (\sig_\Sigma - \sigeps) \zeta \,,
	\label{intquartaeps}
\Esist
for all $V$-valued step functions $z$ and~$\zeta$.
Thus, \accorpa{terzaeps}{quartaeps} hold as well, and the proof of Theorem~\ref{Existenceeps} is complete.


\section{Existence of a Weak Solution}
\label{EXISTENCE}
\setcounter{equation}{0}

In this section, we prove Theorem~\ref{Existence}.
We start from the approximating problem 
\an{analyzed in the previous section}
and let $\eps$ tend to zero.
Since we did not prove uniqueness for the approximating solution,
we take a particular one, namely, the solution we have constructed \an{above}.
This ensures a number of bounds.
Indeed, by the estimates established for the discrete solution and the semicontinuity  of the norms,
it is clear~that
\Bsist
  && \norma\vveps_{\L2\VV}
  + \norma\peps_{\L{4/3}H} 
  + \norma\mueps_{\L2\Vz}
  \non
  \\
  && \quad {}
  + \norma\phieps_{\L\infty V\cap\L2W}
  + \norma\sigeps_{\W{1,4/3}\Vp\cap\L\infty H\cap\L2V}
  \non
  \\
  && \quad {}
  + \eps^{1/2} \, \norma\mueps_{\L\infty H}
  + \eps^{1/2} \, \norma{\dt\phieps}_{\L2H}
  \leq \an{c}\,,
  \label{stimaeps}
\Esist
\an{for a positive constant $c$ independent of $\eps$.}
However, we need \juerg{some additional} estimates.

\step
Sixth a priori estimate

We take any $\badphi\in\L2\Vz$, and (\aat) we test \eqref{secondaeps}, written at the time $t$, by~$\badphi(t)$.
Then, we integrate over~$(0,T)$.
With a procedure that is completely similar to the one used to prove~\eqref{perquintastima},
we see~that
\Beq
  \Bigl| \ioT \< \dt(\eps\mueps+\phieps)(t) , \badphi(t) >_{\Vz} \Bigr|
  \leq c \, \norma\badphi_{\L2\Vz} \, ,
  \non
\Eeq
\an{meaning} that
\Beq
  \norma{\dt(\eps\mueps+\phieps)}_{\L2\Vzp} \leq c \,.
  \label{sestastima}
\Eeq

\step
Seventh a priori estimate

We recall that $\Feps'=\betaeps+\pi$ and test \eqref{terzaeps} written at the time $t$ by $\betaeps(\phieps(t))$.
Then, we integrate over~$(0,T)$.
By also accounting for \last{\eqref{hpzero}, \eqref{disugBetaeps}, \eqref{hpheps} and \eqref{stimaeps}}, we obtain~that
\Bsist
  && \eps \iO \Betaeps(\phieps(T))
  + \intQ \betaeps'(\phieps) |\nabla\phieps|^2
  + \intQ |\betaeps(\phieps)|^2 
  \non
  \\
  && = \eps \iO \Betaeps(\phiz)
  + \intQ \bigl( \mueps \an{{}+ \heps{}} - \Nph(\phieps,\sigeps)  - \pi(\phieps)\bigr) \betaeps(\phieps)
  \non
  \\
  && \leq \frac 12 \intQ |\betaeps(\phieps)|^2 
  + c \,.
  \non
\Esist
Since all \juerg{of} the terms on the \lhs\ are nonnegative, we conclude that
\Beq
  \norma{\betaeps(\phieps)}_{\L2H} \leq c \,.
  \label{settimastima}
\Eeq

\step
Conclusion of the proof of Theorem~\ref{Existence}

By recalling \accorpa{stimaeps}{settimastima}, we see that\an{, as $\eps \to 0$,}
\Bsist
  && \vveps \to \vv
  \quad \hbox{weakly in $\L2\VV\emb\L2{\LLx4}$}\,,
  \label{convvveps}
  \\
  && \peps \to p
  \quad \hbox{weakly in $\L{4/3}H$}\,,
  \label{convpeps}
  \\
  &&\mueps \to \mu
  \quad \hbox{weakly in $\L2\Vz$}\,,
  \label{convmueps}
  \\
  &&\phieps \to \phi
  \quad \hbox{weakly star in $\L\infty V\cap\L2W$}\,,
  \label{convphieps}
  \\
  && \betaeps(\phieps) \to \xi
  \quad \hbox{weakly in $\L2H$}\,,
  \label{convxieps}
  \\
  &&\sigeps \to \sig
  \quad \hbox{weakly star in $\W{1,4/3}\Vp\cap\L\infty H\cap\L2V$}\,,
  \qquad
  \label{convsigeps}
  \\
  && \eps\mueps \to 0 
  \quad \hbox{strongly in $\L\infty H \an{{}\cap \L2{ V_0} {}}$}\,,
  \label{epsmu}
  \\
  && \eps\dt\phieps \to 0
  \quad \hbox{strongly in $\L2H$}\,,
  \label{epsdtphi}
  \\
  &&\eps\mueps+\phieps \to \phi
  \quad \hbox{weakly in $\H1\Vzp\cap\L2V$}\,,
  \label{convcombination}
\Esist
for suitable limit functions \an{$\soluz$}.
More precisely, this holds for some subsequence $\eps_n\searrow0$.
Nevertheless, here and in the sequel, we write \an{just} $\eps$ for simplicity.
Then, the sextuple $\soluz$ satisfies \an{the regularity properties in} \Regsoluz.
By noting that $(\eps\mueps+\phieps)(0)$ converges to $\phi(0)$ weakly in~$\Vzp$,
we infer that the initial conditions~\eqref{cauchy} are satisfied as well,
and we now prove that the sextuple we have found \an{yields in fact a solution to the original system} \Pbl.
To this end, we try to follow the lines used 
to solve the approximating problem at the end of the previous section.
However, the \pier{analogues} of some of those convergence properties require some \an{further} work.

First of all, we clearly have that\an{, as $\eps\to 0$,}
\Beq
  \Nph(\phieps,\sigeps) \to \Nph(\phi,\sigma)
  \aand
  \Ns(\phieps,\sigeps) \to \Ns(\phi,\sigma)
  \quad \hbox{weakly in $\L2V$}.
  \non
\Eeq
Now, we account for the strong compactness results \an{\cite[Thm.~5.1, p.~58]{Lions} and \cite[Sect.~8, Cor.~4]{Simon} to}
deduce that 
$\eps\mueps+\phieps$ converges to $\phi$ strongly in~$\L2H$.
By combining this with \eqref{epsmu}, we infer~that\an{, as $\eps\to 0$,}
\Beq
  \phieps \to \phi
  \quad \hbox{strongly in $\L2H$}.
  \label{strongphieps}
\Eeq
On the other hand, we have the interpolation identity
\Beq
  (\Hx2,\an{H})_{1/8} 
  = \Hx{7/4}
  \non
\Eeq
and the associated interpolation inequality
\Beq
  \norma v_{\Hx{7/4}}
  \leq c \, \norma v_{\Hx2}^{1/8} \, 	\an{\norma v^{7/8}}
  \quad \hbox{for every $v\in\Hx2$} \,.
  \non
\Eeq
Therefore, by applying the \Holder\ inequality, we deduce, for every $v\in\L2{\Hx2}$, that
\Bsist
  && \ioT \norma{v(t)}_{\Hx{7/4}}^2 \, dt
  \leq \ioT \norma{v(t)}_{\Hx2}^{1/4} \, \an{\norma{v(t)}}^{7/4} \, dt
  \non
  \\
  && \leq \Bigl( \ioT \norma{v(t)}_{\Hx2}^2 \, dt \Bigr)^{1/8}
    \Bigl( \ioT \an{\norma{v(t)}}^2 \, dt \Bigr)^{\an{7/8}} \,.
  \non
\Esist
By applying this inequality to $\phieps-\phi$ 
and owing to the boundedness in $\L2W$ and to strong convergence in $\L2H$,
we \an{readily} deduce that
\Beq
  \phieps \to \phi
  \quad \hbox{strongly in $\L2{\Hx{7/4}}$}.
  \non
\Eeq
Now, we recall the (three-dimensional) \an{continuous} embedding
\Beq
  \Hx s \emb \Wx{1,q},
  \quad \hbox{where} \quad 
  s > 1 , \quad
  q \geq 2,
  \aand
  1 - \frac 3q \leq s - \frac 32, 
  \non
\Eeq
and apply it with $s=7/4$ and $q=4$.
We conclude that\an{, as $\eps\searrow 0$,}
\Beq
  \phieps \to \phi
  \quad \hbox{strongly in $\L2{\Wx{1,4}}$}
  \non
\Eeq
\last{which} is the analogue of \eqref{strongphik}.
As in the previous proof, we derive the analogues of \eqref{strongsigk}, \an{\eqref{weakprod:1}--\eqref{weakprod:2}},
and the first of~\eqref{strongnonlin}, namely,
\Bsist
  && \sigeps \to \sigma
  \quad \hbox{strongly in $\L2{\Lx4}$ and \aeQ}\,,
  \non
  \\
  && \an{\nabla\phieps \cdot \vveps \to \nabla\phi \cdot \vv
  \quad \hbox{weakly in $\L1{\an{H}}$}}\,,
  \non
	\\
	&&
	\an{\mueps \nabla\phieps \to \mu \nabla\phi
  \quad \hbox{weakly in $\L1{\an{\HH}}$}}\,,\,
  \non
  \\
  && \an{\Ss(\phieps,\sigeps) \to \Ss(\phi,\sig)}
  \quad \hbox{strongly in $\L2{\Lx4}$}.
  \non
\Esist
Moreover, \an{by} the same \an{argument}, we obtain the analogue of \eqref{convprod},
which now sounds
\Beq
  \intQ (\nabla\sigeps \cdot \vveps) \zeta
  \to \intQ (\nabla\sig \cdot \vv) \zeta 
  \quad \hbox{for every $\zeta\in\L\infty V$}\,.
  \non
\Eeq
\an{In a similar fashion}, dealing with the $\Ns$ terms in place of the velocities, 
one obtains the analogue of \eqref{convprodbis},~i.e.,
\Beq
  \intQ \an{ \Ns(\phieps,\sigeps)\nabla\sigeps \cdot \bzeta
  \to \intQ \Ns(\phi,\sigma)\nabla\sigma \cdot \bzeta}
  \quad \hbox{for every $\bzeta\in\L\infty\VV$} \,.
  \non
\Eeq
At this point, we recall that the approximating solution we are considering also satisfies
\eqref{eqdiveps} and the four variational equations \eqref{intprimaeps} and
\accorpa{intsecondaeps}{intquartaeps}, the first of which being satisfied
for every $\bzeta\in\L\infty\VV$ 
and the others for arbitrary step functions $\an{\badphi}$, $z$ and $\zeta$ 
taking values in $\Vz$, $V$ and~$V$, respectively.
At this point, we can let $\eps$ tend to zero.
By recalling \eqref{hpgeps} and \eqref{hpheps},
we obtain \eqref{eqdiv} and integrated version of the variational equations
\eqref{prima} and \accorpa{seconda}{quarta} with the same \an{type} of test functions.
Thus, those equations are satisfied as they are written in the original problem.
It remains to verify \eqref{xibetaphi}.
To this end, it suffices to observe that 
the weak convergence \eqref{convxieps} coupled with the strong convergence \eqref{strongphieps}
allows us to apply a well-known property of the Yosida approximation
(see, e.g.,  \cite[Prop.~2.2, p.~38]{Barbu})
which yields \an{the inclusion} $\xi\in\beta(\phi)$\an{,} as desired.
Thus, the proof of Theorem~\ref{Existence} is complete.

\Brem
\label{Regphi}
We \an{further} justify \an{the regularity properties claimed in} \eqref{morereg}.
The rigorous argument should involve a regularization of $\beta$
and \an{truncation} in the choice of the test function.
However, we proceed formally, for brevity.
We write \eqref{terza} in the form
\Beq
  \iO \nabla\phi \cdot \nabla z
  + \iO \beta(\phi) z
  = \iO f z
  \quad \hbox{for every $z\in V$} 
  \label{foremorereg}
\Eeq
where $f:=\mu+h-\Nph(\phi,\sig)-\pi(\phi)$.
Then, we test the above equation by $z=(\beta(\phi))^5$ (\aet)
and owe to the Young inequality \eqref{young} on the \rhs.
We obtain (\aet)~that
\Beq
  \iO 5 |\beta(\phi)|^4 |\nabla\phi|^2
  + \iO |\beta(\phi)|^6
  = \iO |f| \,|\beta(\phi)|^5
  \leq \frac 56 \iO |\beta(\phi)|^6
  + \frac 16 \iO |f|^6
  \non
\Eeq
whence immediately $\norma{\beta(\phi)}_6\leq\norma f_6$.
Since $f\in\L2V\subset\L2{\Lx6}$, 
we deduce that $\beta(\phi)\in\L2{\Lx6}$.
Then, elliptic regularity \an{theory} yields $\phi\in\L2{\Wx{2,6}}$.
To show the remaining regularity property,
we test the variational equation \eqref{foremorereg} by $-\Delta\phi$
and recall that $\Nph(\phi,\sig)=-\chi\sigma$.
We thus have~that
\Bsist
  && \norma{\Delta\phi}^2
  + \iO \beta'(\phi) |\nabla\phi|^2
  \leq \norma{\nabla\mu} \, \norma{\nabla\phi} 
  + \norma{h+\chi\sigma-\pi(\phi)} \, \norma{\Delta\phi} 
  \non
  \\
  && \leq \norma{\nabla\mu} \, \norma{\nabla\phi}
  + \frac 12 \, \norma{\Delta\phi}^2
  + c \, \norma{h+\chi\sigma-\pi(\phi)}^2 \,.
  \non
\Esist
Since $\norma{\nabla\phi}$ and the last norm are bounded over~$(0,T)$ 
by \eqref{regh}, \last{\eqref{regphi} and~\eqref{regsigma}}, 
we deduce~that
\Beq
  \frac 12 \, \norma{\Delta\phi}^2
  \leq c \, \norma{\nabla\mu}
  + c
  \quad \hbox{whence also} \quad
  \norma{\Delta\phi}^4
  \leq c \, \norma{\nabla\mu}^2
  + c \,.
  \non
\Eeq
Then, \eqref{regmu} implies \pier{that} $\Delta\phi\in\L4H$,
\pier{\last{and} elliptic regularity theory ensures that} $\phi\in\L4\Hdue$.  
\Erem

\vskip 6mm
\noindent{\an{\bf Acknowledgements}}

\noindent
\an{This research was supported by the Italian Ministry of Education, 
University and Research (MIUR): Dipartimenti di Eccellenza Program (2018--2022) 
-- Dept.~of Mathematics ``F.~Casorati'', University of Pavia. 
In addition, {PC and AS gratefully acknowledge some other support 
from the MIUR-PRIN Grant 2020F3NCPX ``Mathematics for industry 4.0 (Math4I4)'' and}
their affiliation to the GNAMPA (Gruppo Nazionale per l'Analisi Matematica, 
la Probabilit\`a e le loro Applicazioni) of INdAM (Isti\-tuto 
Nazionale di Alta Matematica). }

\End{document}


\an{
\begin{remark} \label{REM:PHEN}
Let us now comment on how the first estimate may be adjusted to overcome also the case when the source terms $\Sph$ and $\Ss$ possess the explicit form in \eqref{phen:laws}.
The only terms that need to be bounded in a different way are the ones examined in the estimate \eqref{est:source}.
Recalling that $P$ is nonnegative and bounded, we substitute estimate \eqref{est:source} with the following 
\begin{align*}
	&  \int_{Q_t} \Sph(\ph,\s,\m) \m + \int_{Q_t} \Ss(\ph,\s,\m) N_\s 
	\\ & \quad 
	 = - \int_{Q_t} P(\ph) \big|\s + \chi(1-\ph)-\m \big|^2
	 -\int_{Q_t} \Sph(\ph,\s,\m) g
	\\ & \quad 
	 \leq 
	-\int_{Q_t} \Sph(\ph,\s,\m) g 
	= - \int_{Q_t} P(\ph)(\s + \chi(1-\ph))g
	+ \int_{Q_t} P(\ph)\m  g
	\\ & \quad  
	\leq 
	 \d \int_{Q_t}|{\nabla\m}|^2 
	 + \cd \int_{Q_t}(|{\ph}|^2+|{g}|^2+|{\s}|^2+1) 
	\\ & \quad  
	 \leq 
	 \d \int_{Q_t}|{\nabla\m}|^2 
	 + \cd \int_{Q_t}(|{\ph}|^2 +|{\s}|^2+1),
\end{align*}
for a positive $\d$ to be fixed at the end of the computations and where we also use the Poincar\'e inequality.
This allows us to consider also the form \eqref{phen:laws} in Theorem \ref{THM:EXWEAK}. The passage to the limit and all the other estimates can be reproduced in the same fashion.
\end{remark}
}


Furthermore, the embeddings $V\emb H$ and $H\emb\Vp$ are compact,
so that we obtain from Ehrling's lemma the compactness inequality 
\Beq
  \norma v
  \,\leq\, \delta\, \norma{\nabla v}
  + C_\delta \, \norma v_{\Vp}
  \quad \hbox{for every $v\in V$ and $\delta>0$},
  \label{compact}
\Eeq
with some $C_\delta>0 $ that depends only on~$\Omega$ and~$\delta$.


as well as the Gagliardo--Nirenberg inequalities we briefly recall.
If $s,q,r\in[1,+\infty]$, $m,j\in\enne$ with $0\leq j < m$, and $\theta \in [\frac jm, 1]$ satisfy the relation
\begin{align}
    j- \frac ds = \left(m -\frac dr\right)\theta + \left(- \frac dq \right)(1-\theta)
    \label{hpGN}
\end{align}
there exists a positive constant $C$ depending only on $\Omega$ and the above parameters such~that
\begin{align}
    \norma{{D^j v}}_s 
	\leq C \, \norma v_{m,r}^\theta \norma v_q^{1-\theta} 
	\quad \text{for all $v\in\Wx{m,r}\cap\Lx q$}.
    \label{GN}
\end{align}

\footnotesize
\begin{thebibliography}{99}

\bibitem{ACHE}
A. Agosti, P. Colli, H. Garcke and E. Rocca.
A Cahn--Hilliard model coupled to viscoelasticity with large deformations.
Preprint arXiv:2204.04951 [math.AP], 1-48, 2022.


\an{
\bibitem{AbelsTer}
H. Abels \pier{and} Y. Terasawa.
On Stokes operators with variable viscosity in bounded and unbounded domains. 
{\it Math. Ann.}, {\bf 344}, 381--429, 2009.}

\bibitem{Barbu}
V. Barbu.
``Nonlinear Differential Equations of Monotone Type in Banach Spaces'',
Springer,
London, New York, 2010.

\bibitem{Brezis}
H. Brezis.
``Op\'erateurs maximaux monotones et semi-groupes de contractions
dans les espaces de Hilbert'',
North-Holland Math. Stud.
{\bf 5},
North-Holland,
Amsterdam,
1973.
 
\bibitem{Byr}
H.M. Byrne, J.R. King, D.L.S. McElwain, and L. Preziosi. 
A two-phase model of solid tumour growth. 
{\it Appl. Math. Lett.} {\bf 16}, 567--573, 2003.


\bibitem{CGH}
P. Colli, G. Gilardi, and D. Hilhorst.
On a Cahn--Hilliard type phase field system related to tumor growth,
{\it Discrete Contin. Dyn. Syst.} {\bf 35}, 2423--2442, 2015.

\pcol{\bibitem{CGRS}
P. Colli, G. Gilardi, E. Rocca, and J. Sprekels.
Vanishing viscosities and error estimate for a Cahn--Hilliard type phase field system related to tumor growth,
{\it Nonlinear Anal. Real World Appl.} \textbf{26}, 93--10, 2015.}


\bibitem{CL}
V. Cristini and J. Lowengrub.  
Multiscale Modeling of Cancer: An Integrated Experimental and Mathematical Modeling Approach. Cambridge University Press, 2010. 
 
\bibitem{EGAR}
M. Ebenbeck and H. Garcke.
Analysis of a Cahn--Hilliard--Brinkman model for tumour growth with chemotaxis,
{\it J. Differential Equations} \textbf{266}, 5998--6036, 2019.

\bibitem{EGAR2}
M. Ebenbeck and H. Garcke.
On a Cahn--Hilliard--Brinkman model for tumor growth and its singular limits.
{\it SIAM J. Math. Anal.} {\bf 51},  1868--1912, 2019.



\bibitem{Fri}
A. Friedman. 
Free boundary problems for systems of Stokes equations. 
{\it Discrete Contin. Dyn. Syst. Ser. B}, 
{\bf 21}, 1455--1468, 2016.

\bibitem{Galdi}
G.P. Galdi.
An Introduction to the Mathematical Theory of the Navier-Stokes
Equations. Springer Monographs in Mathematics, 2nd ed. New York: Springer, Steady-state
problems, 2011.

\bibitem{GL1}
H. Garcke and K.F. Lam.
Well-posedness of a Cahn--Hilliard system modelling tumour growth with chemotaxis and active transport.
{\it  \pier{European J. Appl. Math.}} {\bf 28}, 284--316, 2017.


\bibitem{GARL_2}
H. Garcke and K.F. Lam.
Analysis of a Cahn--Hilliard system with non-zero Dirichlet 
conditions modeling tumor growth with chemotaxis,
{\it Discrete Contin. Dyn. Syst.} {\bf 37}, 4277--4308, 2017.
%
\bibitem{GARL_4}
H. Garcke and K.F. Lam.
On a Cahn--Hilliard--Darcy system for tumour growth 
with solution dependent source terms. 
In ``Trends on Applications of Mathematics to Mechanics'', 
E.~Rocca, U.~Stefanelli, L.~Truskinovski, A.~Visintin~(ed.), 
{\it Springer INdAM Series} {\bf 27}, Springer, Cham, 243--264, 2018.

\bibitem{GLSS}
H. Garcke, K.F. Lam, E. Stika, and V. Styles.
A Cahn--Hilliard--Darcy model for tumour growth with chemotaxis and active transport. 
{\it Math. Models Methods Appl. Sci.} {\bf 26}, 1095--1148, 2016.

 

\bibitem{GLS}
H.~Garcke, K.F.~Lam, and A.~Signori.
On a phase field model of Cahn--Hilliard type for tumour growth with mechanical effects. 
{\it Nonlinear Anal. Real World Appl.}  {\bf 57}, \pier{Paper No. 103192, 28 pp.}, 2021. 

\bibitem{FGR}
S. Frigeri, M. Grasselli, and E. Rocca.
On a diffuse interface model of tumor growth,
{\it  European J. Appl. Math.\/} {\bf 26}, 215--243, 2015. 

\bibitem{FLR}
S. Frigeri, K.F. Lam, and E. Rocca. 
On a diffuse interface model for tumour growth with non-local interactions and degenerate mobilities.
In ``Solvability, Regularity, Optimal Control 
of Boundary Value Problems for PDEs'', 
P.~Colli, A.~Favini, E.~Rocca, G.~Schimperna, J.~Sprekels~(ed.), 
Springer INdAM Series~{\bf 22}, Springer, Cham, 217--254, 2017.

\an{\bibitem{FLS}
S.~Frigeri, K.F.~Lam, and A.~Signori.
Strong well-posedness and inverse identification problem of a non-local phase field tumor model with degenerate mobilities. 
{\it European J. Appl. Math.}, {\bf 33}, 267--308, 2022.}%

\pier{\bibitem{Giga}
Y. Giga.
Analyticity of the semigroup generated by the {S}tokes operator in {$L_{r}$} spaces.
{\it Math. Z.}, {\bf 178}, 297--329, 1981.}

\bibitem{HD}
A. Hawkins-Daarud, K. G. van der Zee, and J. T. Oden. 
Numerical simulation of a thermodynamically consistent four-species tumor growth model.
{\it Int. J. Numer. Math. Biomed. Eng.} {\bf 28}, 3--24, \pier{2012}.  
%
  
\an{
\bibitem{KS2}
P. Knopf and A. Signori. 
Existence of weak solutions to multiphase Cahn--Hilliard--Darcy and Cahn--Hilliard--Brinkman models for stratified tumor growth with chemotaxis and general source terms. 
{\it Comm. Partial Differential Equations}, \pier{{\bf 47}, 233--278, 2022.}}

\juerg{
\bibitem{KRS}
P. Krej\v{c}\'i, E. Rocca, and J. Sprekels. Analysis of a tumor model as a multicomponent deformable porous medium.
{\it Interfaces Free Bound.}, published online first: 2022-03-29, DOI: 10.4171/IFB/472.
}


\bibitem{Lions}
J.-L. Lions.
``Quelques M\'ethodes de R\'esolution des Probl\`emes aux Limites non Lin\'eaires'', 
Dunot, Gauthier-Villars, Paris, 1969. 

\bibitem{LioMag}
J.L. Lions and E. Magenes.
``Non-Homogeneous Boundary Value Problems and Applications'',
Die Grundlehren der Mathematischen Wissenschaft, {\bf 181}, 
Springer-Verlag, Berlin, 1972.



\an{%
\bibitem{RSchS}
E.~Rocca, G.~Schimperna, and A.~Signori.
On a Cahn--Hilliard--Keller--Segel model with generalized logistic source describing tumor growth.
Preprint arXiv:2202.11007 [math.AP], 1--38, 2022.}

\bibitem{SS}
L.~Scarpa and A.~Signori.
On a class of non-local phase-field models for tumor growth with possibly singular potentials, chemotaxis, and active transport.
{\it Nonlinearity}  {\bf 34}, 3199--3250, 2021. 

\bibitem{Simon}
J. Simon.
Compact sets in the space $L^p(0,T;B)$. 
{\it Ann. Mat. Pura Appl.} (4) {\bf 146} 56--96, 1987. 

\bibitem{Tartar}
L. Tartar.
``Remarks on some interpolation spaces",
Carnegie Mellon University,
Research Report No. 94-NA-002,
1994.

\end{thebibliography}
